\documentclass[titlepage,11pt]{article} 
\usepackage{latexsym}
\usepackage{amsfonts}
\usepackage{amsmath}
\newcommand\blackslug{\hbox{\hskip 1pt \vrule width 4pt height 8pt depth 1.5pt
        \hskip 1pt}}
\newcommand\bbox{\hfill \quad \blackslug \bigbreak}
\def\DD{\hbox{-}}
\def\LL{,\ldots,}
\def\cupcup{\cup\cdots\cup}
\newcommand{\vare}{\varepsilon}
\usepackage{tikz}
\usepackage{float}
\usetikzlibrary{decorations.pathreplacing,decorations.markings}

\title{Pure pairs. X. Tournaments and the strong Erd\H{o}s-Hajnal property.\thanks{This is an accepted manuscript. The version of record appeared in European Journal of Combinatorics, Volume 115, January 2024, 103786 at \url{https://doi.org/10.1016/j.ejc.2023.103786}.}}
\author{Maria Chudnovsky\thanks{Supported by  NSF grant DMS 1763817.}\\
Princeton University, Princeton, NJ 08544
\\
\\
Alex Scott\thanks{Research supported by EPSRC grant EP/V007327/1.}\\
Mathematical Institute, University of Oxford, Oxford OX2 6GG, UK
\\
\\
Paul Seymour\thanks{Supported by AFOSR grants
A9550-19-1-0187 and FA9550-22-1-0234, and by NSF grants  DMS-1800053 and DMS-2154169.}\\
Princeton University, Princeton, NJ 08544
\\
\\
Sophie Spirkl\thanks{We acknowledge the support of the Natural Sciences and Engineering Research
Council of Canada (NSERC), [funding reference number RGPIN-2020-03912].
Cette recherche a \'et\'e financ\'ee par le Conseil de recherches en sciences
naturelles et en g\'enie du Canada (CRSNG), [num\'ero de r\'ef\'erence
RGPIN-2020-03912].  }\\
University of Waterloo, Waterloo, Ontario N2L3G1, Canada}

\date{March 13, 2021; revised \today}

\newtheorem{thm}{}[section]

\newcommand{\Proof}{\noindent{\bf Proof.}\ \ }

\begin{document}
\maketitle
\begin{abstract} 
A {\em pure pair} in a tournament $G$ is an ordered pair $(A,B)$ of disjoint subsets of $V(G)$ such that every vertex in $B$ is adjacent from 
every vertex in $A$. Which tournaments $H$
have the property that if $G$ is a tournament not containing $H$ as a subtournament, and $|G|>1$, there is a pure pair $(A,B)$
in $G$ with $|A|,|B|\ge c|G|$, where $c>0$ is a constant independent of $G$? Let us say that such a tournament $H$ has the
{\em strong EH-property}

As far as we know, it might be that a tournament $H$ has this property if and only if its vertex set has a linear ordering in which
its backedges form a forest. Certainly this condition is necessary, but we are far
from proving sufficiency. We make a small step in this direction, showing that if a tournament can be ordered with at most three backedges
then it has the strong EH-property (except for one case, that we could not decide). In particular, every tournament with at most 
six vertices has the property, except for three
that we could not decide. We also give a seven-vertex 
tournament that does not have the strong EH-property.

This is related to the Erd\H{o}s-Hajnal conjecture, which in one form says that for 
every tournament $H$ there exists $\tau>0$
such that
every tournament $G$ not containing $H$ as a subtournament has a transitive subtournament of cardinality at least $|G|^\tau$.
Let us say that a tournament $H$ satisfying this has the {\em EH-property}. It is known that every tournament with the strong 
EH-property also has the EH-property; so our result extends work
by Berger, Choromanski and Chudnovsky, who proved that every tournament with at most six vertices has the EH-property, except for one
that they did not decide.
\end{abstract}

\section{Introduction}

A {\em tournament} is a digraph $G$, with no loops, such that for every pair $u,v$ of distinct vertices, exactly one of $uv,vu$
is an edge. (All graphs and digraphs in this paper are finite, and have no loops or parallel edges.) Let us say a tournament $G$ {\em contains} a tournament $H$
if there is a subtournament of $G$ isomorphic to $H$, and $G$ is {\em $H$-free} otherwise. We denote the number of vertices of $G$
by $|G|$.

The {\em Erd\H{o}s-Hajnal conjecture} was raised as a question by 
Erd\H{o}s and Hajnal~\cite{EH0,EH} and asserts that, for every graph $H$, there exists $\tau>0$ such that every graph $G$ not containing
an induced subgraph isomorphic to $H$ has a clique or stable set of cardinality at least $|G|^\tau$.  Alon, Pach and
Solymosi~\cite{aps} showed that it is equivalent to the following assertion about tournaments:
\begin{thm}\label{apsconj}
{\bf Conjecture:} For every tournament $H$ there exists $\tau>0$ such that every $H$-free tournament $G$ has a transitive subtournament
with at least $|G|^\tau$ vertices.
\end{thm}
For a tournament $H$, if there exists $\tau>0$ as in \ref{apsconj}, we say that $H$ has the {\em EH-property} or {\em (weak) EH-property}.
Thus the conjecture says that all tournaments have the EH-property.

Let $P_7$ denote the Paley tournament with seven vertices; that is, its vertex set is $\{1\LL 7\}$, and for all distinct
$i,j\in \{1\LL 7\}$, $j$ is adjacent from $i$ if $j-i$ is congruent to $1,2$ or $4$ modulo 7. Let $P_7^-$
denote the tournament obtained by deleting one 
vertex from $P_7$. (It makes no difference which vertex is deleted.)
Berger, Choromanski and Chudnovsky~\cite{berger} showed:
\begin{thm}\label{EH6}
Every tournament with at most six vertices has the EH-property, 
except possibly for $P_7^-$.
\end{thm}
There are other classes of tournaments that have been shown to have the EH-property: see for instance~\cite{galaxies, constellations, 7vertex}.

In a graph $G$, a {\em pure pair} is a pair $A,B$ of disjoint subsets of $V(G)$ such that either there are no edges between $A,B$
or every vertex in $A$ is adjacent to every vertex in $B$; and its {\em order} is $\min(|A|,|B|)$.
Let us say a {\em pure pair} in a tournament $G$ is an ordered pair $(A,B)$ of disjoint subsets of $V(G)$ such that every vertex in $B$
is adjacent from every vertex in $A$; and its {\em order} is $\min(|A|,|B|)$. And let us say a tournament $H$ has the 
{\em strong EH-property} or {\em SEH property} if there exists $c>0$
such that for every $H$-free tournament $G$ with $|G|>1$, there is a pure pair in $G$ with order at least $c|G|$.
It is easy to see that every tournament with the strong EH-property also has the EH-property, but not
all tournaments have the strong EH-property; we shall see that $P_7$ does not.
So it is natural to ask which tournaments do. (One can take the same approach for graphs -- see~\cite{pure1}.) 
This question seems not to have been studied to any great extent. We discussed it briefly in~\cite{pure1}; ``heroes'', 
defined in~\cite{heroes}, have the strong EH-property; and 
Berger, Choromanski, Chudnovsky and Zerbib~\cite{bergerstrong} proved that $D_5$ has the strong EH-property ($D_5$ is defined below);
but we know of nothing else on the topic.

A {\em numbering} of a graph is an enumeration $(v_1\LL v_n)$ of its vertex set; and an  {\em ordered graph} is a graph
together with some numbering.
If $(v_1\LL v_n)$ is a numbering of a tournament $H$, then the corresponding {\em backedge graph} 
of $H$ is the ordered graph $B$ with vertex set $V(H)$ and numbering $(v_1\LL v_n)$, 
in which for $1\le i<j\le n$, $v_i$ and $v_j$ are adjacent 
in $B$ if and only if $v_i$ is adjacent from $v_j$ in $H$. Its edges are called {\em backedges}.
A tournament can be reconstructed from a backedge graph
and the corresponding numbering, and it is often convenient to work with the backedge graph rather than directly with the
tournament.

Different numberings of the same tournament may result in wildly different
backedge graphs, of course.
For instance, $D_5$ is the (unique, up to isomorphism) tournament with five vertices, in which every vertex has outdegree two;
and the following are two of its backedge graphs (in such figures, vertices are always 
numbered from left to right):
\begin{figure}[H]
\centering

\begin{tikzpicture}[scale=1,auto=left]

\tikzstyle{every node}=[inner sep=1.5pt, fill=black,circle,draw]

\node (v1) at (1,0) {};
\node (v2) at (2,0) {};
\node (v3) at (3,0) {};
\node (v4) at (4,0) {};
\node (v5) at (5,0) {};
\draw (v1) to [bend right=20] (v5);
\draw (v1) to [bend left=20] (v4);
\draw (v2) to [bend left=20] (v5);

\node (u1) at (1,-1) {};
\node (u2) at (2,-1) {};
\node (u3) at (3,-1) {};
\node (u4) at (4,-1) {};
\node (u5) at (5,-1) {};
\draw (u1) to [bend left=20] (u3);
\draw (u1) to [bend left=20] (u5);
\draw (u3) to [bend left=20] (u5);
\draw (u2) to [bend right=20] (u4);

\end{tikzpicture}

\caption{Two backedge graphs for $D_5$.} \label{fig:D_5}
\end{figure}
For a graph $G$, we denote its complement graph by $\overline{G}$; and if $G$ is an ordered graph, $\overline{G}$ means the complement graph with the same numbering.
Let us say the {\em reverse} $\overline{H}$ of a tournament
$H$ is obtained by reversing the direction of all edges of $H$. A tournament has the strong EH-property if and only if
its reverse does. 
Note that, if under some numbering a tournament $H$ has backedge graph $B$, then $\overline{B}$ is the backedge graph of 
$\overline{H}$ under the same numbering; and $B$ with its numbering reversed is 
the backedge graph of $\overline{H}$
under the reverse numbering. 

Can we hope to characterize the tournaments with the strong EH-property? A parallel question for graphs had a very satisfactory answer:
we proved in~\cite{pure1} that:
\begin{thm}\label{pure1}
For a graph $H$, the following are equivalent:
\begin{itemize}
\item there exists $c>0$ such that for every graph $G$ with $|G|>1$ not containing $H$ or $\overline{H}$ as an induced subgraph, 
there is a pure pair $A,B$ in $G$ with order at least $c|G|$;
\item one of $H,\overline{H}$ is a forest.
\end{itemize}
\end{thm}
One might hope for a parallel for this in the world of tournaments. Certainly, one half is true:
we will show, in \ref{forest2}, that
\begin{thm}\label{forest}
Every tournament with the strong EH-property admits a numbering for which the backedge graph is a forest.
\end{thm}
As far as we know, the converse to this might also be true. Initially this seemed unlikely to us, but we have tried hard to disprove 
it and failed, so let us pose it as a conjecture:
\begin{thm}\label{forestconj}
{\bf Conjecture: }A tournament has the strong EH-property if and only if it admits a numbering for which the backedge graph is a forest.
\end{thm}
This would be a beautiful analogue of \ref{pure1}, but we are far from
proving it. Indeed, the tournament $P_7^-$ has a backedge graph that is a forest with only four edges (see figure \ref{fig:P7-}), and we cannot even show
that it has the (weak) EH-property.

\begin{figure}[H]
\centering

\begin{tikzpicture}[scale=1,auto=left]

\tikzstyle{every node}=[inner sep=1.5pt, fill=black,circle,draw]

\node (v1) at (1,0) {};
\node (v2) at (2,0) {};
\node (v3) at (3,0) {};
\node (v4) at (4,0) {};
\node (v5) at (5,0) {};
\node (v6) at (6,0) {};
\draw (v1) to [bend right=20] (v6);
\draw (v1) to [bend left=20] (v4);
\draw (v3) to [bend left=20] (v6);
\draw (v2) to [bend right=20] (v5);

\end{tikzpicture}

\caption{Backedge graph for $P_7^-$.} \label{fig:P7-}
\end{figure}
There is a recent positive result in this area: it is shown in~\cite{density4} that if a tournament can be built from nothing by
repeatedly adding vertices with in-degree at most one or out-degree at most one (and consequently admits an ordering in which the
backedge graph is a forest), then it has the (weak) EH-property.

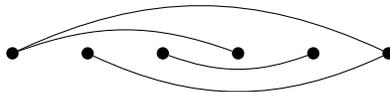
\begin{figure}[H]
\centering

\begin{tikzpicture}[scale=1,auto=left]

\tikzstyle{every node}=[inner sep=1.5pt, fill=black,circle,draw]

\node (v1) at (1,0) {};
\node (v2) at (2,0) {};
\node (v3) at (3,0) {};
\node (v4) at (4,0) {};
\node (v5) at (5,0) {};
\node (v6) at (6,0) {};
\draw (v1) to [bend left=25] (v6);
\draw (v1) to [bend left=20] (v4);
\draw (v2) to [bend right=25] (v6);
\draw (v3) to [bend right=20] (v5);

\end{tikzpicture}

\caption{Backedge graph for $H_6$.} \label{fig:H6}
\end{figure}

The first main result of this paper is:
\begin{thm}\label{threeedges}
Let $H$ be a tournament that admits a numbering $(v_1\LL v_n)$ for which the backedge graph $B$ has at most three edges.
Suppose that $H$ is also $D_5$-free, that is, there do not exist $a,b,c,d, e$ with $1\le a<b<c<d<e\le n$ such that $E(B)=\{v_av_d, v_av_e, v_bv_e\}$.
Then $H$ has the strong EH-property.
\end{thm}
To see the equivalence asserted in the second sentence, observe that if there exist $a,b,c,d,e$ as stated then $H$ contains $D_5$;
and conversely, if $H$ contains $D_5$ then there exist $a,b,c,d,e$ as stated, since
up to isomorphism there is only one backedge graph of $D_5$ with only three edges. (We leave the reader to check this.)
Perhaps the second sentence in \ref{threeedges} (the condition about $D_5$) can be omitted, but that remains open.
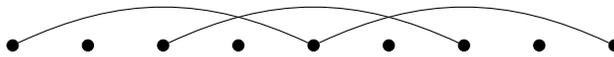
\begin{figure}[H]
\centering

\begin{tikzpicture}[scale=1,auto=left]

\tikzstyle{every node}=[inner sep=1.5pt, fill=black,circle,draw]

\node (v1) at (1,0) {};
\node (v2) at (2,0) {};
\node (v3) at (3,0) {};
\node (v4) at (4,0) {};
\node (v5) at (5,0) {};
\node (v6) at (6,0) {};
\node (v7) at (7,0) {};
\node (v8) at (8,0) {};
\node (v9) at (9,0) {};
\draw (v1) to [bend left=25] (v5);
\draw (v3) to [bend left=25] (v7);
\draw (v5) to [bend left=25] (v9);

\end{tikzpicture}

\caption{Backedge graph of a tournament satisfying \ref{threeedges}.} \label{fig:Refex}
\end{figure}
For instance, \ref{threeedges} implies that the tournament with the backedge graph in figure \ref{fig:Refex} has the strong 
EH-property. A referee kindly told us that the methods of earlier papers would not show this.
Let $H_6$ be the tournament with a backedge graph as in figure \ref{fig:H6}. 

\noindent \ref{threeedges} will be used to show our second main result, that:
\begin{thm}\label{mainthm}
Every tournament with at most six vertices has the strong EH-property, except possibly for $P_7^-, H_6$ and $\overline{H_6}$.
\end{thm}
All except one of them (and that one is easy) either contain $D_5$ or admit backedge graphs with at most
three edges, and so we can apply \ref{threeedges}.

We remark that if we just wanted to prove that these tournaments have the (weak) EH-property, we could make use of the theorem of 
Alon, Pach and Solymosi~\cite{aps} that the class of tournaments with the EH-property is closed under substitution, and so
it would only be necessary to examine the tournaments that are not built from smaller ones by substitution.
But the class with the strong EH-property is {\em not} closed under substitution.

The paper is organized as follows. Sections \ref{sec:SEH}--\ref{sec:threeback} are devoted to proving \ref{threeedges}, and then we 
turn to \ref{mainthm}. We need to prove that if $|H|\le 6$ then $H$ has the SEH property (except for three cases).
\begin{itemize}
\item  In sections \ref{sec:sparsity} and \ref{sec:6vertexwithD5}, we prove that if $H$ contains $D_5$ then $H$ has the SEH property (except for two of the exceptional cases).
This proof is
a modification of a proof of Berger, Choromanski, Chudnovsky and Zerbib~\cite{bergerstrong}, who showed that $D_5$ itself has the
strong EH-property.
\item In section \ref{sec:6vertexnoD5}, we prove by case-by-case analysis, that every tournament with at most six vertices 
admits a numbering 
for which the backedge graph has at most three edges, except for four particular tournaments, three of which are the exceptions
in \ref{mainthm} (it is easy to show that the fourth has the SEH property). So \ref{mainthm} follows from
\ref{threeedges}.
\item In section \ref{sec:forests} we prove \ref{forest}, and show that $P_7$ does not
have the SEH property. Finally, in section \ref{sec:P7bad} we give two other tournaments that have the SEH property but not a ``rainbow''
refinement of it discussed in the proof of \ref{threeedges}.
\end{itemize}

\section{The strong EH-property for ordered graphs}\label{sec:SEH}

An ordered graph $G$ {\em contains} another ($H$ say)
if some induced subgraph of $G$, with the induced numbering, is isomorphic (as an ordered graph) to $H$; and if not, $G$ is {\em $H$-free}.
It is helpful to recast our problem about tournaments into the language of ordered graphs. 
We observe first that:
\begin{thm}\label{purepair}
Let $G$ be a tournament, and let $J$ be a backedge graph of $G$, with numbering $(v_1\LL v_n)$.
\begin{itemize}
\item If $G$ has a pure pair of order $t$, then $J$ has a pure pair of order at least $t/2$.
\item If $J$ has a pure pair of order $t$ then $G$ has a pure pair of order at least $t/2$.
\end{itemize}
\end{thm}
\Proof
Let $(A,B)$ be a pure pair of order $t$ in $G$, and choose $i\in \{1\LL n\}$ minimum such that one of $A',B'$
has cardinality at least $t/2$, where $A'=\{v_1\LL v_i\}\cap A$ and $B'=\{v_1\LL v_i\}\cap B$. Define $A''=A\setminus A'$
and $B''=B\setminus B'$.
If $|A'|\ge t/2$, then from the minimality of $i$,
$|B'|< t/2$, and so $|B''|\ge t/2$;
and since every vertex in $B$ is $G$-adjacent from every vertex in $A$, it follows that there are no edges of $J$
between $A'$ and $B''$, and this pair of sets is the desired pure pair of $J$.
Similarly, if $|B'|\ge t/2$, then every vertex of $B'$ is $J$-adjacent to every vertex in
in $A''$, and so this is the desired pure pair. This proves the first assertion.

For the second, let $A,B$ be a pure pair of order $t$ in $J$, choose $i$ as before, and define $A',A'', B', B''$ as before.
By exchanging $A,B$ if necessary,
we may assume that $|A'|\ge t/2$; and so either $(A',B'')$ (if there are no edges
of $J$ between $A,B$) or $(B'',A')$ (if every vertex in $A$ is $J$-adjacent to every vertex in $B$) is the
desired pure pair of $G$. This proves the second assertion, and so proves \ref{purepair}.~\bbox

Alon, Pach and Solymosi~\cite{aps} proved that the Erd\H{o}s-Hajnal conjecture is equivalent to the same statement for ordered
graphs: that is, for every ordered graph $H$ there exists $\tau>0$ such that every $H$-free ordered graph $G$
has a clique or stable set of cardinality at least $|G|^\tau$. One can extend the
``strong EH-property'' to ordered graphs in the natural way, but while it makes sense to ask which tournaments have the strong
EH-property, the same question for ordered graphs is unprofitable, as only very trivial ordered graphs have the property. For instance,
 a result of Fox~\cite{fox} shows that the ordered graph with vertices $v_1,v_2,v_3$ numbered in this order, and edges $v_1v_2,v_2v_3$, does {\em not} have the property
(see \cite{pure5, pure6} for related results).
We can show (we omit the proof) that if an ordered graph has this property, then each of its components either has at most two
vertices, or is a three-vertex path with middle vertex the first or last in the induced numbering, or is
one particular four-vertex ordered path.

It is better to exclude more than one ordered graph at the same time.
Let $\mathcal{A}$  be a set of ordered graphs. We say an ordered graph $G$ is {\em $\mathcal{A}$-free} if $G$
is $H$-free for all $H\in \mathcal{A}$; and $\mathcal{A}$ has the {\em strong EH-property} if there exists $c>0$
such that every $\mathcal{A}$-free ordered graph $G$ with $|G|>1$ has a pure pair of order at least $c|G|$.

To translate our question about tournaments into the language of ordered graphs, we observe that
\begin{itemize}
\item because of \ref{purepair}, a tournament $G$ has a linear pure pair if and only if some (or equivalently, every)  backedge graph of $G$
has a linear pure pair (with a different constant of linearity);
\item
a tournament $G$ does not contain a tournament $H$
if and only if some (and therefore every) backedge graph of $G$ contains none of $B_1\LL B_k$, where $B_1\LL B_k$ are the backedge
graphs of $H$ that arise from the different numberings of $H$.
\end{itemize}
Thus a tournament $H$
has the strong EH-property if and only if the set $\mathcal{A}$ of all backedge graphs that arise from $H$ under its different numberings
has the strong EH-property.

This set $\mathcal{A}$ can be rather large. For instance, when $H$ is $D_5$, it has 24 nonisomorphic backedge graphs, and we are looking at the ordered graphs that
contain none of 24 specific ordered graphs. Excluding just one of them is not enough, but 24 is more than we need; the proof
given in \ref{D5mainthm} shows that 
a subset of four of them already has the strong EH-property, the two
shown in figure \ref{fig:D_5} and their complements.
A similar thing happens for all the tournaments we can handle: we need to retain at most three (usually two) backedge graphs and their complements.

\section{Blockades and rainbow subgraphs}

Our goal at the moment is to show that all $D_5$-free tournaments that admit backedge graphs with at most
three edges have the strong EH-property. We will prove that in fact they have a stronger property that we explain now.

A {\em blockade} in a set $V$ is a family $\mathcal{B}=(B_i: i\in I)$ of pairwise disjoint nonempty subsets of $V$, where $I$ is a
finite set of integers. (We have used blockades in several papers of this series, for instance in~\cite{pure1}.)
Its {\em length} is $|I|$,
and the minimum of $|B_i|\;(i\in I)$ is its {\em width}. We write $W(\mathcal{B})$ to denote the width of $\mathcal{B}$.
We call the sets $B_i$ {\em blocks} of the blockade.
(What matters is that the blocks are not too small. We could shrink the larger ones to make them all the same size.)
We are interested in blockades of some fixed length in the vertex set of some graph, ordered graph or tournament $G$, 
in which each block contains
linearly many vertices of $G$.

If $(v_1\LL v_n)$ is a numbering of $V$, a blockade $(B_i: i\in I)$ {\em respects} the numbering if
for all $i_1,i_2\in I$ with $i_1<i_2$, if $v_h\in B_{i_1}$ and $v_j\in B_{i_2}$
then $h<j$. In this case we say $\mathcal{B}$ is {\em respectful}.

Let $\mathcal{B}$ be a blockade in a set $V$. A graph (or ordered graph, or tournament) $H$ with $V(H)\subseteq V$ is
{\em $\mathcal{B}$-rainbow} 
if each vertex of $H$ belongs to some block of $\mathcal{B}$, and no two vertices belong to the same block.
A {\em copy} of a graph (or ordered graph, or tournament) $H$ is another such object isomorphic to $H$.

In order to prove that a tournament has the strong EH-property, it is often easier to prove something even stronger.
Let us say 
 a tournament $H$ has the {\em rainbow strong EH-property} or {\em RSEH-property} if there exists $c$ with $0<c<1$ 
such that if $\mathcal{B}$ is a blockade of length at least $1/c$ in a tournament $G$, and there is no $\mathcal{B}$-rainbow 
copy of $H$ contained in $G$,
then there is a pure pair in $G$ of order at least $cW(\mathcal{B})$.
\begin{thm}\label{rainbowstrongtour}
If $H$ is a tournament with the RSEH-property then $H$ has the strong EH-property.
\end{thm}
\Proof
Choose $c>0$ as in the definition of the  RSEH-property; by reducing $c$ we may assume that $k=1/c$ is an integer. 
Let $c'=c^2/2$. Now let $G$ be an $H$-free tournament with $|G|>1$.
We claim that $G$ has a pure pair of order at least $c'|G|$. Since $|G|>1$, we may assume that $|G|>1/c'$, since otherwise
a pure pair of order 1 exists and satisfies the theorem. There is a blockade $\mathcal{B}$ in $G$
of length $k$, where
$$W(\mathcal{B})\ge \lfloor|G|/k\rfloor=\lfloor c|G|\rfloor\ge c|G|/2;$$ 
and since $G$ is $H$-free there is certainly no
$\mathcal{B}$-rainbow copy of $H$ contained in $G$. Thus $G$ has a pure pair of order at least $cW(\mathcal{B})\ge c^2|G|/2=c'|G|$.
This proves \ref{rainbowstrongtour}.~\bbox

The converse of \ref{rainbowstrongtour} is not true: we will see that $D_5$ has the strong EH-property, but not the RSEH-property.
On the other hand, we will show that:
\begin{thm}\label{3edgerainbow}
If $H$ is a $D_5$-free tournament that admits a backedge graph with at most three edges, then 
$H$ has the RSEH-property.
\end{thm}

Similarly, if $\mathcal{A}$ is a set of ordered graphs, we say that $\mathcal{A}$ has the  {\em rainbow strong EH-property} or 
{\em RSEH-property} if there exists $c$ with $0<c<1$
such that if $\mathcal{B}$ is a respectful blockade of length at least $1/c$ in an ordered graph $G$, 
and there is no $\mathcal{B}$-rainbow  
copy of any member of $\mathcal{A}$ contained in $G$,
then there is a pure pair in $G$ of order at least $cW(\mathcal{B})$.

Evidently we have:
\begin{thm}\label{ordertotour}
If $H$ is a tournament and the set of all backedge graphs of $H$ (or a subset of this set) has the  RSEH-property, 
then $H$ has the  RSEH-property.
\end{thm}

An {\em anticomplete pair} in a graph $G$ (possibly ordered) is a pair $A,B$ of disjoint subsets of $V(G)$ such that there are no edges
between $A,B$; and its {\em order} is $\min(|A|,|B|)$.

We need to throw a form of sparsity into this sea of definitions too: we say that a set $\mathcal{A}$ of ordered graphs
has the {\em sparse rainbow strong EH-property} or {\em SRSEH-property} if there exists $c$ with $0<c<1$
such that if $\mathcal{B}$ is a respectful blockade of length at least $1/c$ in an ordered graph $G$, 
and there is no $\mathcal{B}$-rainbow    
copy of any member of $\mathcal{A}$ contained in $G$,
and every vertex of $G$ has degree less than $cW(\mathcal{B})$, 
then there is an anticomplete pair in $G$ of order at least $cW(\mathcal{B})$. We say $c$ is an {\em SRSEH-coefficient} for $\mathcal{A}$.

Next we show that if a set of ordered graphs has the SRSEH-property then it together with its set of complement graphs
has the RSEH-property. The proof will 
use the following theorem of \cite{EHC5}:
\begin{thm}\label{foxrodl}
For all $\vare>0$ and every graph $P$ on $p$ vertices, there exist $\gamma, \delta>0$ such that
if $G$ is a graph containing
fewer than $\gamma|G|^p$ induced labelled copies of $P$, then there exists $X\subseteq V(G)$ with $|X|\ge \delta|G|$
such that one of $G[X], \overline{G}[X]$ has maximum degree at most $\vare\delta|G|$.
\end{thm}

\begin{thm}\label{SRSEH}
Let $\mathcal{A}$ be a set of ordered graphs, and let $\mathcal{A}'$ be the set of complements of the members of $\mathcal{A}$.
If $\mathcal{A}$ has the SRSEH-property then $\mathcal{A}\cup \mathcal{A}'$ has the RSEH-property.
\end{thm}
\Proof
Choose $H\in \mathcal{A}$. 
By \ref{orderedornot}, there is a graph $P$ such that for
every numbering of $P$, the ordered graph that results contains $H$.
Let $p=|P|$.

Let $c'$ be an SRSEH-coefficient for  $\mathcal{A}$.
By reducing $c'$ we may assume that $1/c'$ is an integer at least two ($1/c'=K'$ say).

Let 
$\vare\le c'^2/2$ with $\vare>0$,
and choose $\gamma, \delta>0$ to satisfy \ref{foxrodl}. Choose $c$ with $1/c$ an integer ($1/c=K$ say),
such that $c\le c'\delta/2$, and 
$(1-cp)^p>1-\gamma$,
and
$c\le \delta/c-1/c'$. We claim that $c$ satisfies our requirement.
\\
\\
(1) {$\delta K/K'-1\ge \max\left(\vare \delta K/c', c/c'\right)$.}
\\
\\
To see that $\delta K/K'-1 \ge \vare \delta K/c'$, observe that $\delta K/(2K')=\delta c'/(2c)\ge 1$, and
$\delta K/(2K')\ge \vare \delta K/c'$. The second part, that $\delta K/K'-1 \ge c/c'$, is true from the choice of $c$.
This proves (1).

\bigskip

Let $G$ be an ordered graph with a blockade $\mathcal{B}=(B_1\LL B_K)$ that respects the numbering $(v_1\LL v_n)$ of $G$, such that
there is no $\mathcal{B}$-rainbow copy of any member of $\mathcal{A}\cup \mathcal{A}'$. Let 
$W=W(\mathcal{B})$. We must show that 
there is a pure pair in $G$ of order at least $cW$.
We may assume that $V(G)=B_1\cupcup B_K$, and that 
$|B_i|=W$ for each $i$, and so $|G|=KW$, and 
$$B_i=\{v_j:\;(i-1)W<j\le iW\}$$ 
for $1\le i\le K$.

The number of injections $\phi$ from $V(P)$ into $V(G)$ such that the vertices $\phi(v)\;(v\in V(P))$
all belong to different blocks of the blockade, is
$$|G|(|G|-W)(|G|-2W)\cdots (|G|-(p-1)W)> (1-p/K)^p|G|^p\ge (1-\gamma)|G|^p,$$
and since none of them give an isomorphism from $P$ to an induced subgraph of $G$ (from the choice of $P$, and since 
there is no $\mathcal{B}$-rainbow copy of $H$ in $G$),
it follows that the number of induced labelled copies of $P$
in $G$ is less than $|G|^p-(1-\gamma)|G|^p=\gamma|G|^p$.
By \ref{foxrodl}, there exists $X\subseteq V(G)$ with $|X|\ge \delta KW$, such that
one of $G[X], \overline{G}[X]$ has maximum degree less than $\vare\delta KW$; and by replacing $G$ by $\overline{G}$ if
necessary, we may assume that $G[X]$ has maximum degree less than $\vare\delta KW$.
By (1), there exists a real number $W'$ such that
$$\frac{\delta K}{K'}-1\ge \frac{W'}{W}\ge \max\left(\vare \delta K/c',\frac{c}{c'}\right).$$

The sets $B_1\cap X\LL B_K\cap X$ each have cardinality at most $W$, but their union has cardinality at least $\delta KW$.
Define $i_0=0$, and inductively for $s=1,2,\ldots$ choose $i_s\in \{1\LL K\}$ minimum such that
$|B_s'|\ge W'$, where $B_s'=\bigcup_{i_{s-1}<i\le i_s}B_i\cap X$, if such a choice is possible; and let the first value of $s$
where the choice is impossible be $s=t+1$. Thus $B_1'\LL B_t'$ are defined.
From the minimality of each $i_s$ it follows that $|B_s'|\le W'+W$
for $1\le s\le t$. From the maximality of $t$,
$$|X\cap \left(B_{i_t+1}\cupcup B_K\right)|<W';$$
and so $\delta KW\le |X|\le t(W'+W)+W'$. Since $\delta K/K'-1\ge W'/W$ it follows that $t\ge K'$.

Let $\mathcal{B}'$ be the blockade $(B_1'\LL B_{K'}')$; it has width at least $W'$, and it respects the numbering $(v_1\LL v_n)$.
Let $G'=G[B_1'\cupcup B_{K'}']$. Every vertex of $G'$ has degree less than $\vare\delta KW\le c'W'$ in $G'$.
Also there is no $\mathcal{B}'$-rainbow copy of any member of $\mathcal{A}$ in $G'$, since such a copy would also be 
$\mathcal{B}$-rainbow. Hence from the choice of $c'$, there is an anticomplete pair in $G'$ of order at least $c'W'\ge cW$.
This proves \ref{SRSEH}.~\bbox


By combining \ref{SRSEH}, \ref{ordertotour} and \ref{rainbowstrongtour}, we have:
\begin{thm}\label{summary}
Let $H$ be a tournament. If $\mathcal{A}$ is a set of some of the backedge graphs of $H$, and $\mathcal{A}$ has the SRSEH-property,
then $H$ has the RSEH-property and hence the strong EH-property.
\end{thm}
\Proof
Let $\mathcal{A}'$ be the set of complements of the members of $\mathcal{A}$. Thus all the members of $\mathcal{A}'$ are also backedge graphs of $H$ under appropriate numberings of $H$ (obtained by reversing the numberings that give the members of $\mathcal{A}$).
Since $\mathcal{A}$ has the SRSEH-property, \ref{SRSEH} implies that $\mathcal{A}\cup \mathcal{A}'$ has the RSEH-property; and hence
so does the set of all backedge graphs of $H$. Consequently $H$ has the RSEH-property, by \ref{ordertotour}, and hence the strong EH-property, by 
\ref{rainbowstrongtour}. This proves \ref{summary}.~\bbox

\section{Blockades}

If we start with a blockade with great length, we might hope to make a smaller, but more tightly structured, blockade by 
shrinking or removing some
of its blocks.
Here are two useful ways to make smaller blockades from larger.
First, if $\mathcal{B}=(B_i:i\in I)$ is a blockade, let $I'\subseteq I$; then $(B_i:i\in I')$ is a blockade, of smaller length
but of at least the same width, and we call it a {\em sub-blockade} of $\mathcal{B}$.
Second, for each $i\in I$ let $B_{i}'\subseteq B_{i}$ be nonempty; then
the sequence $(B_i':i\in I)$ is a blockade, of the same length but possibly of smaller width, and we call it a {\em contraction}
of $\mathcal{B}$. A contraction of a sub-blockade (or equivalently, a sub-blockade of a contraction) we call a {\em minor} of $\mathcal{B}$.


If $H$ is a $\mathcal{B}$-rainbow induced subgraph,
its {\em support} is the set of all $i\in I$ such that $V(H)\cap B_i\ne \emptyset$.
If $J$ is an ordered graph, we define the {\em trace} of $J$ (relative to $\mathcal{B}=(B_i:i\in I)$) to be the set of supports of all $\mathcal{B}$-rainbow copies of $J$.
If $\tau\ge 1$ an integer, we say $\mathcal{B}$ is {\em $\tau$-support-uniform} if
for every ordered graph $J$ with $|J|\le \tau$, either the
trace of $J$ is empty, or it consists of all subsets of $I$ of cardinality $|J|$.

Let $0< \kappa\le 1$ and $\tau\ge 1$.
We say $\mathcal{B}$ is
{\em $(\kappa,\tau)$-support-invariant} if
for every contraction $\mathcal{B}'=(B_i':i\in I)$ of $\mathcal{B}$ of width at least $\kappa$ times the width of $\mathcal{B}$,
and for every ordered graph $J$ with $|J|\le \tau$,
the trace of $J$ relative to $\mathcal{B}$ equals the trace of $J$ relative to $\mathcal{B'}$.

We need a theorem of \cite{pure1}, the following:
\begin{thm}\label{getminor}
Let $k\ge 0$ and $\tau\ge 1$ be integers, and $0<\kappa\le 1$; then there exist an integer $K$ with the following property.
Let $\mathcal{B}=(B_1\LL B_K)$ be a blockade in a graph.
Then there is a minor $\mathcal{B}'$ of $\mathcal{B}$, with
length $k$ and width at least $\kappa^{2^{K+\tau^2}}W(\mathcal{B})$, such that
$\mathcal{B}'$ is $\tau$-support-uniform and $(\kappa,\tau)$-support-invariant.
\end{thm}

In fact we have cheated a little here: the theorem of \cite{pure1} defines  ``$\tau$-support-uniform'' and
``$(\kappa,\tau)$-support-invariant'' using ordered trees rather than general ordered graphs. But it is not worth writing the proof out again, since exactly the same argument works
for general ordered graphs, except we have to replace the bound $\tau^\tau$ used in \cite{pure1} for the number of
ordered trees
on at most $\tau$ vertices, by the bound $2^{\tau^2}$ for the number of ordered graphs on at most $\tau$ vertices. (This
only changes the multiplicative constant in \ref{getminor}.)

One important application of \ref{getminor} is the following.
\begin{thm}\label{separate}
Let $\mathcal{A}$ be a set of ordered graphs, and let $H\in \mathcal{A}$. Suppose that $V_1,V_2\subseteq V(H)$ 
with $V_1\cup V_2=V(H)$ and $V_1\cap V_2=\emptyset$, such that there are no edges of $H$ between $V_1,V_2$. For $i = 1,2$,
let $H_i$ be the ordered subgraph of $H$ induced on $V_i$, and let
$\mathcal{A}_i=\{H_i\}\cup (\mathcal{A}\setminus \{H\})$. If both $\mathcal{A}_1,\mathcal{A}_2$ have the SRSEH-property then so does
$\mathcal{A}$.
\end{thm}
\Proof
For $i = 1,2$, let $c_i$ be an SRSEH-coefficient for $\mathcal{A}_i$.
By reducing
$c_1$ or $c_2$ we may assume that $c_1=c_2$ ($=c'$ say), and $1/c'\ge 2|H|$, and $1/c'$ is an integer. Let $k=1/c'$, and choose $K$ as in \ref{getminor},
taking $\tau=|H|$ and $\kappa=1/2$.  Let $c>0$ with $c\le c'2^{-2^{K+|H|^2}}$ and $c\le 1/K$. 

We claim that,
if $\mathcal{B}$ is a respectful blockade of length at least $1/c$ in an ordered graph $G$, and there is no $\mathcal{B}$-rainbow
copy in $G$ of any member of $\mathcal{A}$, and every vertex of $G$ has degree less than $cW(\mathcal{B})$,
then there is a pure pair in $G$ of order at least $cW(\mathcal{B})$. Let $W=W(\mathcal{B})$.
By \ref{getminor}, there is a minor $\mathcal{C}=(C_1\LL C_k)$ of $\mathcal{B}$, of width $W'$ where $W'\ge 2^{-2^{K+|H|^2}}W$,
such that $\mathcal{C}$ is $|H|$-support-uniform and $(1/2,|H|)$-support-invariant. If there is no $\mathcal{C}$-rainbow
copy in $G$ of any member of $\mathcal{A}_1$, then since $G$ has maximum degree less than $cW\le c'W'$, there is an anticomplete
pair in $G$ of order at least $c'W'\ge cW$ as required. So we may assume that there is a $\mathcal{C}$-rainbow copy in $G$ of some member
of $\mathcal{A}_1$. But there is no $\mathcal{B}$-rainbow
copy in $G$ of any member of $\mathcal{A}$, so there is a $\mathcal{C}$-rainbow copy in $G$ of $H_1$, and similarly
we may assume (for a contradiction) that there is a  $\mathcal{C}$-rainbow copy in $G$ of $H_2$. 

Let $(v_1\LL v_n)$ be the numbering of $H$, and let $I_j= \{i : v_i \in V_j\}$ for $j=1,2$.
Since  $\mathcal{C}$ is 
$|H|$-support-uniform, there is a $\mathcal{C}$-rainbow copy $J_1$ of $H_1$ with support $I_1$.
For $1\le i\le k$, if $i\notin I_2$ let $D_i=C_i$, and if $i\in I_2$ 
let $D_i$ be the set of vertices in $C_i$ with no neighbour in
$V(J_1)$. Thus $|D_i|\ge |C_i|-cW|J_1|\ge |C_i|/2$, since
$$cW|J_1|\le cW|H|\le (c'|H|)\left(2^{-2^{K+|H|^2}}W\right)\le W'/2\le |C_i|/2.$$
Let $\mathcal{D}=(D_1\LL D_k)$.
Since $\mathcal{C}$ is $|H|$-support-uniform and $(1/2,|H|)$-support-invariant, and there is a $\mathcal{C}$-rainbow 
copy of $H_2$ in $G$, it follows
that there is a $\mathcal{D}$-rainbow copy $J_2$ of $H_2$ in $G$ with support $I_2$. But there are no edges of $G$
between $V(J_1)$ and $V(J_2)$, and so $G$ contains a $\mathcal{C}$-rainbow copy of $H$, a contradiction. This proves
\ref{separate}.~\bbox

Let $\mathcal{A}$ be a set of ordered graphs, and for each $H\in \mathcal{A}$ let $C_H$ be a component of $H$ 
(with the induced numbering). We call the set $\{C_H:H\in \mathcal{A}\}$ a {\em transversal} of $\mathcal{A}$.
By repeated application of \ref{separate}, it follows that:

\begin{thm}\label{connected}
Let $\mathcal{A}$ be a finite set of ordered graphs. If every transversal of $\mathcal{A}$ has the SRSEH-property then
$\mathcal{A}$ has the SRSEH-property.
\end{thm}

For instance, let $\mathcal{A}=\{H_1,H_2\}$, where for $i = 1,2$, $H_i$ is an ordered graph in which two components
($A_i,B_i$ say) have more than two vertices, and perhaps some other components have at most two vertices.
In order to prove that $\mathcal{A}$ has the SRSEH-property it would suffice to show
that 
$$\{A_1,A_2\},\{A_1,B_2\},\{B_1,A_2\},\{B_1,B_2\}$$
all have the SRSEH-property. There are other transversals, using a one- or two-vertex component of one of $H_1$ or $H_2$, but they all 
obviously have the SRSEH-property and we don't have to check them. But in general we do have to check all four of the transversals 
given.
When we come to work with a tournament, the art will be to select a small number of orderings of the tournament 
such that {\em every} transversal
of the corresponding set $\mathcal{A}$ of backedge graphs has the SRSEH-property, so that we can apply \ref{connected}.

\section{Some sets of ordered graphs that have the SRSEH-property}

In this section we prove that certain sets of ordered graphs have the SRSEH-property, enough that for every tournament 
we need to handle, it has a set of backedge graphs such that all their transversals can be shown to have the SRSEH-property.

A {\em left-star} is an ordered graph, with numbering $(v_1\LL v_n)$ say, such that $n> 0$, and $v_1$ is adjacent to every
other vertex, and every edge is incident with $v_1$. If it has $n$ vertices it is also called a {\em left $(n-1)$-star}.
We start with an easy one:
\begin{thm}\label{cliques}
Let $\mathcal{A}$ be a set of ordered graphs, such that all components of some member of $\mathcal{A}$ are left-stars or right-stars,
and all components of some member of $\mathcal{A}$ are cliques. Then $\mathcal{A}$ has the SRSEH-property.
\end{thm}
\Proof
By \ref{connected} it suffices to prove the result when $\mathcal{A}$ has two members, one a left-star and one a clique. Choose
$t$ such that both these ordered graphs have at most $t$ vertices.

Let $N\ge 0$ be an integer such that every graph on at least $N$ vertices has either a stable set or a clique of size $t-1$.
Let $c=1/(N+1)$. Let $G$ be an ordered graph, let $\mathcal{B}=(B_1\LL B_K)$ be a respectful blockade
in $G$ of length at least $1/c$, and suppose that there is no $\mathcal{B}$-rainbow copy in $G$
of any member of $\mathcal{A}$, and every vertex of $G$ has degree less than $cW(\mathcal{B})$. (The last condition will not be used.)
We claim that $G$ has an anticomplete pair of order at least $cW(\mathcal{B})$. We may assume that each $B_i$ has cardinality $W(\mathcal{B})$.
For $2\le i\le N+1$, we may assume that fewer than $cW(\mathcal{B})$ vertices in $B_1$ have no neighbour in $B_i$, since this set of vertices
forms an anticomplete pair with $B_i$. Since $Nc<1$, it follows that there is a vertex $v_1\in B_1$ with a neighbour $v_i\in B_i$ for $2\le i\le N+1$.
From the choice of $N$, either $t-1$ of the vertices $v_2\LL v_{N+1}$ form a stable set (and then $G$ contains a $\mathcal{B}$-rainbow 
copy of the left-star in $\mathcal{A}$) or $t-1$ of them form a clique (and then $G$ contains a $\mathcal{B}$-rainbow copy of 
the clique in $\mathcal{A}$), in either case a contradiction. This proves 
\ref{cliques}.~\bbox

We need another theorem of~\cite{pure1}, the following:
\begin{thm}\label{rainbow}
For every tree $T$, there exists $c>0$, such that
for every graph $G$ with a blockade $\mathcal{B}$ of length at least $1/c$, if
there is no $\mathcal{B}$-rainbow copy of $T$, and 
every vertex has degree less than $cW(\mathcal{B})$, 
then 
there is an anticomplete pair of order at least $cW(\mathcal{B})$.
\end{thm}
This is a theorem about unordered graphs, and in particular, the vertices of the $\mathcal{B}$-rainbow copy of $T$ might 
be in any order. Still, we can deduce some useful results about ordered graphs from it, for instance:

\begin{thm}\label{ordertree}
Let $\mathcal{A}$ be a set of ordered graphs that does not have the SRSEH-property. For every tree $T$, there is a 
numbering of $T$ such that the ordered graph formed by $T$ with this numbering contains no member of $\mathcal{A}$.
\end{thm}
\Proof Let $T$ be a tree, and let $c>0$ satisfy \ref{rainbow}. Since $\mathcal{A}$ does not have the SRSEH-property, there is
an ordered graph $G$, and a respectful blockade $\mathcal{B}$ in $G$, such that $G$
has maximum degree less than $cW(\mathcal{B})$, and there is no $\mathcal{B}$-rainbow copy of any member of $\mathcal{A}$ in $G$,
and there is no anticomplete pair in $G$ of order at least $cW(\mathcal{B})$. By the choice of $c$, there is a $\mathcal{B}$-rainbow
copy of the unordered graph $T$ in $G$; let $J$ be this copy, with the induced numbering. Then $J$ contains no
member of $\mathcal{A}$.  This proves \ref{ordertree}.~\bbox

A {\em left-spike}
is an ordered graph, with numbering $(v_1\LL v_n)$ say, such that $n> 1$, and $v_2$ is adjacent to every other vertex, 
and every edge is incident with $v_2$. {\em Right-spikes} 
are defined similarly.
A {\em monotone path}
is an ordered graph,  with numbering $(v_1\LL v_n)$
where $n>0$, with edge set $\{v_1v_2, v_2v_3\LL v_{n-1}v_n\}$.
Left-stars, left-spikes and monotone paths are all special cases of a {\em left-broom}, which is an ordered tree, with numbering $(v_1\LL v_n)$ where $n>0$,
and with edge set
$$\{v_1v_2, v_2v_3\LL v_{m-1}v_m\}\cup \{v_mv_{m+1}, v_mv_{m+2}\LL v_mv_n\}$$
for some $m$ with $1\le m\le n$. A {\em right-broom} is defined similarly.

If $(v_1\LL v_n)$ is a numbering, we say $v_i$ is {\em earlier} than $v_j$ if $i<j$,
and $v_i$ is {\em later} than $v_j$ if $i>j$. Sometimes we will have different graphs and digraphs with the same vertex set, and we will sometimes speak of ``$H$-neighbour'' meaning
``neighbour in $H$'', and so on.
\ref{rainbow} has the following consequence:

\begin{thm}\label{stars}
Let $\mathcal{A}$ be a set of ordered graphs. If 
either $\mathcal{A}$ contains a left-star and a right-broom, or it contains a right-star and a left-broom,
then $\mathcal{A}$ has the SRSEH-property.
\end{thm}
\Proof
Let $\mathcal{A}$ contain a left-star $L$ and a right-broom $R$, both with at most $t\ge 2$ vertices.
Let $T$ be the tree in which every vertex has degree either $2t$ or one; and there is a vertex $v$ with degree $2t$ such that every  
path of $T$ with one end $v$, and maximal with this property, has exactly $t-1$ edges. Take a numbering of $T$, making an ordered graph $T'$. 
For each vertex $u$ of $T'$ with degree $2t$, 
we may assume that at least
$t$ of its neighbours are earlier than $u$, since otherwise $T'$ contains $L$. Let $u_0=v$, 
and having chosen $u_i$, if $i<t$ let $u_{i+1}$ be a $T'$-neighbour of $u_i$ that is earlier than $u_i$.
For each $i$ with $0\le i<t$, there are at least $t$ $T'$-neighbours of $u_i$ that are earlier than $u_i$, and so $T'$
contains $R$ (because $|R|\le t$.) From \ref{ordertree}, this
proves \ref{stars}.~\bbox

This will suffice to handle almost all the tournaments of interest to us, but there are a few tough ones that need
something extra, provided by the following three results.

Let us say a {\em left-bristle} is an ordered graph, with numbering $(v_1\LL v_n)$
where $n>2$, where $v_1$ is adjacent to $v_2\LL v_{n-1}$, and $v_n$ is adjacent to exactly one of $v_2\LL v_{n-1}$,
and there are no other edges. A {\em right-bristle} is defined similarly.
\begin{thm}\label{bristle}
Let $\mathcal{A}$ be a set of ordered graphs containing a left 2-star and a right-bristle; or containing a right 2-star and a left-bristle.
Then $\mathcal{A}$ has the SRSEH-property.
\end{thm}
\Proof
Choose $t\ge 3$ such that $\mathcal{A}$ contains a left 2-star and a right-bristle $R$, both  with at most $t$ vertices.
Let $T$ be the tree in which every vertex has degree either $2t$ or one; and there is a vertex $v$ such that every 
path of $T$ with one end $v$, and maximal with this property, has exactly two edges. Thus $|T|=4t^2+1$. Let $c'$ satisfy \ref{rainbow} (with $c$ replaced by $c'$).

Choose $c>0$ such that $c\le c'/(1+c')$ and $1/c\ge 4t^2+2$. Let $G$ be an ordered graph, let $\mathcal{B}=(B_1\LL B_K)$ be a respectful blockade
in $G$ of length at least $1/c$, and suppose that there is no $\mathcal{B}$-rainbow copy in $G$
of any member of $\mathcal{A}$, and every vertex of $G$ has degree less than $cW(\mathcal{B})$. Let $W=W(\mathcal{B})$.
We claim that $G$ has an anticomplete pair of order at least $cW$. We may assume that each $B_i$ has cardinality $W$.

For $2\le i\le 4t^2+2$ let $C_i$ be the set of vertices in $B_i$ with a neighbour in $B_1$. We may assume that $|B_i\setminus C_i|<cW$,
since $B_i\setminus C_i$ is anticomplete to $B_1$, and so $|C_i|> (1-c)W$. 
If there is no $(C_2\LL C_{4t^2+2})$-rainbow copy of $T$, then by \ref{rainbow}, either some vertex of $G$ has degree at least $c'(1-c)W$,
or $G$ has an anticomplete pair of order at least $c'(1-c)W$, and since $c'(1-c)\ge c$, we are done.

Thus we may assume (for a contradiction) that there is a $(C_2\LL C_{4t^2+2})$-rainbow copy of $T$, and we assume (to simplify notation) that 
this is $T$ itself. Now $v$ has $2t$ neighbours in $T$, and at least $2t-1$ of them are earlier than $v$; let $2t-1$ of them
be $a_1\LL a_{2t-1}$, numbered in order. One neighbour of $a_t$ is later than $a_t$, namely $v$, and so all the others are earlier than
$a_t$, say $b_1\LL b_{2t-1}$, again numbered in order. Since $T$ is $(C_2\LL C_{4t^2+2})$-rainbow, it follows that $b_t$ has a neighbour 
$x\in B_1$. Now $x$ is nonadjacent to $b_1\LL b_{t-1}$ and to $b_{t+1}\LL b_{2t-1}$. If $x$ is also nonadjacent to $a_t$, then the 
ordered subgraph induced on $\{x,a_t, b_1\LL b_{2t-1}\}$ contains $R$, a contradiction. Thus $x$ is adjacent to $a_t$. By the same argument,
$x$ is nonadjacent to $a_1\LL a_{t-1}$ and to $a_{t+1}\LL a_{2t-1}$; and it is also nonadjacent to $v$, since otherwise $\{x,b_t,v\}$
induces a $\mathcal{B}$-rainbow left 2-star. But then $\{x,v, a_1\LL a_{2t-1}\}$ contains $R$, a contradiction. This proves \ref{bristle}.~\bbox

Let us say a {\em crossed left-star} is an ordered graph, with a numbering $(v_1\LL v_n)$, such that $v_1$ is adjacent 
to $v_2\LL v_n$, and only one edge is not incident with $v_1$. (Thus it consists of a left-star and one more 
edge joining some pair of leaves of the left-star.) A {\em crossed right-star} is defined similarly.

\begin{thm}\label{startri}
If $\mathcal{A}$ contains a three-vertex monotone path, and a crossed left-star, and a crossed right-star 
then $\mathcal{A}$ has the SRSEH-property.
\end{thm}
\Proof
Let the crossed left-star and the crossed right-star both have at most $n$ vertices.
Let $c=1/(2n)$, let $G$ be an ordered graph, let $\mathcal{B}=(B_1\LL B_K)$ be a respectful blockade
in $G$ of length $K\ge 1/c$, let $W$ be its width, and suppose that there is no $\mathcal{B}$-rainbow copy in $G$
of any member of $\mathcal{A}$, and every vertex of $G$ has degree less than $cW$. 
We claim that $G$ has an anticomplete pair of order at least $cW$. We may assume that each $B_i$ has cardinality $W$.

For $i=2,3\LL 2n-1$ in turn, we will inductively define $A_i\subseteq B_i$ with the following properties:
\begin{itemize}
\item either at least $(1-1/n)W$ vertices in $B_1$, or at least $(1-1/n)W$
vertices in $B_{2n}$, have a neighbour in $A_i$;
\item the sets $A_2,A_3\LL A_{i}$ are pairwise anticomplete; and
\item for all $j\in \{2\LL 2n-1\}\setminus \{i\}$, fewer than
$cW$ vertices in $B_j$ have a neighbour in $A_i$.
\end{itemize}
The inductive definition is as follows. Assume that $2\le i\le 2n-1$, and $A_2\LL A_{i-1}$ are defined.
Let $X$ be the set of vertices in $B_i$ that have a neighbour in one of $A_{2}\cupcup A_{i-1}$; thus $|X|\le (i-2)cW$.
Since $|B_i\setminus X|\ge cW$, we may assume that fewer than $cW$ vertices in $B_n$ have no neighbour in $B_i\setminus X$, 
since otherwise $G$ has an anticomplete pair of order at least $cW$. Since $cW\le W/n$, we may
choose $A_i\subseteq B_{i}\setminus X$ minimal such that either at least $(1-1/n)W$ vertices in $B_1$, or at least $(1-1/n)$
vertices in $B_{2n}$, have a neighbour in $A_i$. 
Since each vertex in $A_i$ has fewer than  $cW$ neighbours in $B_n$, the minimality of $A_i$ implies that at most
$(1-1/n)W+cW$ vertices in $B_1$ have a neighbour in $A_i$, and hence the set $Z$ of vertices in $B_1$ with no neighbour in $A_i$
has cardinality at least $(1/n-c)W\ge cW$. Similarly the set $Z'$ of vertices in $B_{2n}$ with no neighbour in $A_i$
has cardinality at least $cW$. Let $2\le j\le 2n-1$ with $j\ne i$. We claim that fewer than $cW$ vertices in $B_j$
have a neighbour in $A_i$. To see this, suppose that $j>i$ (the argument when $j<i$ is similar and we omit it).
If $v\in B_j$ has a neighbour in $A_i$,
then it has no neighbour in $Z'$, since there is no $\mathcal{B}$-rainbow monotone three-vertex path in $G$; and since
we may assume that $G$ has no anticomplete pair of order at least $cW$, and $|Z'|\ge cW$, it follows that fewer than $cW$
vertices in $B_j$ have a neighbour in $A_i$ as claimed. This completes the inductive definition. In summary, we have:
\begin{itemize}
\item for $2\le i\le 2n-1$, either at least $(1-1/n)W$ vertices in $B_1$, or at least $(1-1/n)$
vertices in $B_{2n}$, have a neighbour in $A_i$;
\item the sets $A_2,A_3\LL A_{2n-1}$ are pairwise anticomplete; and
\item for all distinct $i,j\in \{2\LL 2n-1\}$, fewer than
$cW$ vertices in $B_j$ have a neighbour in $A_i$.
\end{itemize}
From the symmetry, we may assume that for
at least $n-1$ values of $i\in \{2\LL 2n-1\}$, at least $(1-1/n)W$ vertices in $B_1$
have a neighbour in $A_i$.
Choose $n-1$ such values, say $i_2\LL i_n$ in increasing order, and define $i_1=1$.
Since at most $W/n$ vertices in $B_1$ have no neighbour in each $B_{i_s}$ for $2\le s\le n$, there is a set $M_1\subseteq B_1$ with 
$|M_1|\ge W/n$ such that every vertex in $M_1$ has a neighbour in each of $A_{i_2}\LL A_{i_n}$.

There is a crossed left-star $L$ in $\mathcal{A}$; let its numbering be $(v_1\LL v_t)$ say where $t\le n$, and $v_a,v_b$ are adjacent
for some $a,b$ with $2\le a<b\le t$. Let $M_{i_a}$ be the set of all vertices in $B_{i_a}$ with no neighbour in any of the sets
$A_{i_s}$ where $s\in \{2\LL t\}\setminus \{a,b\}$, and define $M_{i_b}\subseteq B_{i_b}$ similarly. Thus $|M_{i_a}|,|M_{i_b}|\ge W-(t-1)cW\ge 2cW$.
We may assume that fewer than $cW$ vertices in $M_{i_a}$ have no neighbour in $M_1$ (since $|M_1|\ge cW$) and 
similarly, fewer than $cW$ vertices in $M_{i_a}$ have no neighbour in $M_{i_b}$; and so some vertex $u_{i_a}\in B_{i_a}$ 
has a neighbour $u_1\in M_1$
and a neighbour $u_{i_b}\in M_{i_b}$. Since there is no $\mathcal{B}$-rainbow three-vertex path in $G$ it follows that $u_1,u_{i_b}$ 
are adjacent. For each $s\in \{2\LL t\}\setminus \{a,b\}$, choose $u_{i_s}\in A_{i_s}$ adjacent to $u_1$ (this is possible since every vertex of $M_1$
has a neighbour in $A_{i_s}$). Then the ordered subgraph induced on $\{u_{i_1}, u_{i_2}\LL u_{i_t}\}$ is a $\mathcal{B}$-rainbow 
copy of $L$, a contradiction. This proves \ref{startri}.~\bbox

A {\em left-split}  is an ordered graph, with numbering $(v_1\LL v_n)$, such that:
\begin{itemize}
\item $v_1,v_2$ are nonadjacent, and $\{v_3\LL v_n\}$ is a clique; and
\item for $3\le i\le n$, $v_i$ is adjacent to at most one of $v_1,v_2$.
\end{itemize}
\begin{thm}\label{split}
If $\mathcal{A}$ contains a left 2-star, and a crossed right-star, and a left-split,
then $\mathcal{A}$ has the SRSEH-property.
\end{thm}
\Proof
For $t\ge 1$, a {\em $t$-uniform crossed right-star} is a crossed right-star with numbering $(v_1\LL v_{n})$, where $n=3t-1$, and
$v_{t}, v_{2t}$ are adjacent. Every crossed right-star is contained in a $t$-uniform crossed right-star for all sufficiently large $t$.

For $t\ge 1$, a {\em $t$-uniform left-split} is a left-split with numbering $(v_1\LL v_n)$, where $n=3t+2$, such that for $i =3\LL n$, $v_1$
is adjacent to $v_i$ if $i$ is divisible by three, and $v_2$ is adjacent to $v_i$ if $i-1$ is divisible by three.
Every  left-split  is contained in a $t$-uniform left-split for all sufficiently large $t$. Choose $t\ge 1$
such that some member of $\mathcal{A}$ is contained in a $t$-uniform crossed right-star, and some member of $\mathcal{A}$
is contained in a $t$-uniform left-split.

Choose $K$ satisfying \ref{getminor}, taking $k=35t+1$ and $\tau=8t+1$ and $\kappa=1/4$. Let $L=K+(8t+1)^2$.
By \ref{stars},
the set consisting of a left 2-star and a monotone $(8t+1)$-vertex path has the SRSEH-property. Let $c_0>0$ be an SRSEH-coefficient 
for this set. 

Choose $c>0$ with 
$$c\le \min\left(c_04^{-2^{L}},  1/K, 1/(35t+1), 4^{-2^{L}-1}/t\right).$$ 
We will show that $c$ is an SRSEH-coefficient for $\mathcal{A}$.
Let $G$ be an ordered graph, let $\mathcal{B}'$ be a respectful blockade
in $G$ of length at least $1/c$, let $W'$ be its width, and suppose that there is no $\mathcal{B}'$-rainbow copy in $G$
of any member of $\mathcal{A}$. 
We claim that either some vertex has degree at least $cW'$ in $G$, or 
$G$ has an anticomplete pair of order at least $cW'$. 

From the choice of $K$, since $1/c\ge K$, there is a minor $\mathcal{B}=(B_1\LL B_{35t+1})$ of $\mathcal{B}'$ of width at least
$4^{-2^{L}}W'$, such that 
$\mathcal{B}$ is $(8t+1)$-support-uniform and $(1/4,8t+1)$-support-invariant. 
Let its width be $W$; we may assume that all its blocks have cardinality $W$.
There is no $\mathcal{B}$-rainbow left $2$-star, and if there is no $\mathcal{B}$-rainbow monotone $(8t+1)$-vertex path,
then from the choice of $c_0$, there is either some vertex with degree at least $c_04^{-2^{L}}W'$ in $G$, or
an anticomplete pair in $G$ of order at least $c_02^{-2^{L}}W'$; and since $c_04^{-2^{L}}\ge c$, 
in either case this proves our claim.

So we may assume (for a contradiction) that there is a $\mathcal{B}$-rainbow monotone $(8t+1)$-vertex path.
For $i=0,1,2$, let $C_{8ti+1}$ be the set of vertices in $B_{8ti+1}$ that have at least one neighbour in each of 
$B_{32t+2} \LL B_{35t+1}$. For each $j>8ti+1$, at most $cW'$ vertices in $B_{8ti+1}$ have no neighbour in $B_j$; and so $|B_{8ti+1}\setminus C_{8ti+1}|\le 3tcW' $.
Hence $|C_{8ti+1}|\ge W-3tcW'$. Again, for $i=0,1,2$ let $D_i$ be the set of vertices in $B_{24t+1}$ that belong to a 
$\mathcal{B}$-rainbow 
monotone $(8t+1)$-vertex path with support $\{8ti+1\LL 8t(i+1), 24t+1\}$ and with its first vertex (the vertex in $B_{8ti+1}$)
in $C_{8ti+1}$. Now
$|C_{8ti+1}|\ge W-3tcW'\ge W/4$, because $tc\le 4^{-2^{L}-1}\le W/(4W')$.
Since $\mathcal{B}$ is $(8t+1)$-support-uniform and $(1/4,8t+1)$-support-invariant,
and there is no $\mathcal{B}$-rainbow 
monotone $(8t+1)$-vertex path with support $\{8ti+1\LL 8t(i+1), 24t+1\}$ and with its first vertex 
in $C_{8ti+1}$ and last vertex in $B_{24t+1}\setminus D_i$, it follows that $|B_{24t+1}\setminus D_i|<W/4$, and so 
$|D_i|> 3W/4$. 

Let $D_{3}$ be the set of vertices in $B_{24t+1}$ that belong to a $\mathcal{B}$-rainbow monotone $(8t+1)$-vertex path with support
$\{24t+1\LL 32t+1\}$; then by the same argument $|D_{3}|> 3W/4$. Since $D_0\LL D_3$ all have cardinality more than
$3W/4$, and they are all subsets of $B_{24t+1}$ (which has cardinality $W$), there exists $w\in D_0\cap D_1\cap D_2\cap D_3$.
For $i = 0,1,2$, let $P_i$ be a $\mathcal{B}$-rainbow
monotone $(8t+1)$-vertex path with support $\{8ti+1\LL 8t(i+1), 24t+1\}$, with first vertex ($u_i$ say) in $C_{8ti+1}$ and last vertex $w$;
and let $P_3$ be a $\mathcal{B}$-rainbow monotone $(8t+1)$-vertex path with support
$\{24t+1\LL 32t+1\}$ and first vertex $w$.
\\
\\
(1) {\em No vertex in any of $B_{32t+2} \LL B_{35t+1}$ is adjacent to more than one of $u_0,u_1,u_2$.}
\\
\\
Suppose that $z\in B_{32t+2}\cupcup B_{35t+1}$ is adjacent to $u_a, u_b$ say, where $0\le a<b\le 2$. Since there is no $\mathcal{B}$-rainbow
left 2-star, it follows that $z$ is adjacent to every vertex of $P_a\cup P_b$, and hence to every vertex of $P_3$. Let $w'$ be the neighbour
of $w$ in $P_a$.
Since $P_a$ is induced, it has a stable set $I_a$ of cardinality $t$ containing $w'$. Each vertex in $I_a\setminus $ has at most two neighbours in 
$V(P_b)$ since there is no $\mathcal{B}$-rainbow left 2-star; and so there is a stable subset $I_b$ of $P_b$ containing $w$ and with 
cardinality $t$, such that $ww'$ is the only edge of $G[I_a\cup I_b]$. Each vertex in $I_a\cup I_b\setminus \{w\}$ has at most two neighbours in $V(P_3)$;
so there is a stable subset $I_3$ of $P_3$ of cardinality $t$, containing $w$, and such that $ww'$ is the only edge of 
$G[I_a\cup I_b\cup I_3]$. But then
the ordered graph induced on $I_a\cup I_b\cup I_3\cup \{z\}$ is a $\mathcal{B}$-rainbow copy of a $t$-uniform crossed right-star, 
a contradiction. This proves (1).

\bigskip

For $i=3\LL 3t+2$ choose $u_i\in B_{32t-1+i}$, adjacent to $u_1$ if $i$ is divisible by three, adjacent to $u_2$ if $i-1$ is 
divisible by three, and adjacent to $u_0$ otherwise. (This is possible since $u_j\in C_{8tj+1}$ for $j = 0,1,2$.) Thus each of 
$u_3\LL u_{3t+2}$
has exactly one neighbour in $\{u_0,u_1,u_2\}$, by (1). For $3\le i\le 3t+2$, $u_i$ is adjacent to one of $u_0,u_1,u_2$,
and hence to all of one of $V(P_0), V(P_1), V(P_2)$, since there is no $\mathcal{B}$-rainbow left 2-star; and in particular, it is adjacent to $w$.
Thus $w$ is adjacent to each of $u_3\LL u_{3t+2}$, and consequently $\{u_3\LL u_{3t+2}\}$ is a clique; and the ordered subgraph 
induced on $\{u_1\LL u_{3t+2}\}$ is a $t$-uniform left-split, a contradiction.
This proves \ref{split}.~\bbox

\section{Tournaments that have backedge graphs with at most three edges}\label{sec:threeback}

Let us (at last!) apply all these results to prove:
\begin{thm}\label{noD5}
If $H$ is a $D_5$-free tournament with a backedge graph with at most three edges, then $H$ has the RSEH-property.
\end{thm}
\Proof
Let $B$ be a backedge graph of $H$ that has at most three edges, and let $(v_1\LL v_n)$ be its numbering. 
$B$ will have at most six vertices of positive degree, but between them there may
be arbitrary sequences of vertices of degree zero, and we cannot ignore them, because we need to use $B$ to find other backedge graphs 
in order to apply \ref{connected}.  
We need some notation to encompass this. Let $B$ have $t$ vertices of positive degree, and let us number them $b_1,b_3,b_5\LL b_{2t-1}$ in
order; and let us label the sequences of vertices between them as $B_0, B_2\LL B_{2t}$, where 
the sequence $(v_1\LL v_n)$ is the concatenation of 
$$B_0, b_1, B_2, b_3\LL b_{2t-1},B_{2t}.$$
This notation does not tell us the number of vertices in each sequence $B_i$, but we do not need that. 

By \ref{summary} it suffices to show that some set of backedge graphs of $H$ has the SRSEH-property.
First, if no vertex of $B$ has degree more than one, then every component of $B$ is both a left-star and a right-star, and so from
\ref{stars} and \ref{connected}, $\{B\}$ has the SRSEH-property. Thus we may assume that some vertex has degree more than one.
If some $b_i$ is incident with every edge of $B$, then by moving $b_i$ to the start of the numbering, we obtain a numbering
with back-edge graph a left-star (and isolated vertices), and similarly by moving $b_i$ to the other end of the numbering, we obtain
a numbering with back-edge graph a right-star (and isolated vertices), and \ref{stars} and \ref{connected}
imply that the set of these two backedge graphs has the SRSEH-property, and so $H$ has the RSEH-property. So we may assume there is 
no such $b_i$. In particular, $B$ has exactly three edges.

Suppose that $t=3$, and so $b_1,b_3,b_5$ are pairwise adjacent.  
The numberings
$$B_0, b_5, b_1, B_2, b_3,B_4, B_6$$
$$B_0, B_2, b_3,B_4, b_5, b_1, B_6$$
have backedge graphs in which each component is a left-star, and each component is a right-star, respectively. Thus 
each transversal (of this set of two ordered graphs) consists of a left-star and a right-star, and therefore
we may apply \ref{stars} and \ref{connected}.

So we may assume that $t\ge 4$. Suppose that $t=4$. Since $B$ has three edges and no vertex is incident with all of them, 
the subgraph induced on $\{b_1,b_3,b_5, b_7\}$
is a four-vertex path. There are several possibilities for the order in which $b_1,b_3,b_5,b_7$ appear in this path, but there is 
some symmetry we can use to reduce the number of cases. First, there are two orders in which $b_1,b_3,b_5,b_7$ appear in this path, 
reverses of one another, and we only need list one of them. Second, we do not need to list both of two cases which can be taken 
one to the other
by reversing the numbering of $H$; since reversing the numbering of $H$ gives a backedge graph of the reverse of $H$, and 
the result holds for $H$ if and only if it holds for the reverse of $H$.
Up to these two symmetries, the possibilities for the vertices of this path in order are the following:
\begin{itemize}
\item $b_1\DD b_3\DD b_5\DD b_7$. Apply \ref{stars} and \ref{connected} to 
$$B_0, b_1, B_2, B_4, b_5, b_3, B_6, b_7, B_8$$
$$B_0, b_1, B_2, b_5, b_3, B_4, B_6, b_7, B_8.$$
\item $b_1\DD b_3\DD b_7\DD b_5$. Apply \ref{stars} and \ref{connected} to 
$$B_0, b_1, B_2, b_7, b_3, B_4, b_5, B_6, B_8$$
$$B_0, B_2, b_3, b_1, B_4, b_5, B_6, b_7, B_8.$$
\item $b_1\DD b_5\DD b_3\DD b_7$. Apply \ref{stars} and \ref{connected} to
$$B_0, b_5, b_1, B_2, b_3, B_4, B_6, b_7, B_8$$
$$B_0, b_1, B_2, B_4, b_5, B_6, b_7, b_3, B_8.$$
\item $b_1\DD b_5\DD b_7\DD b_3$. Apply \ref{stars} and \ref{connected} to
$$B_0, b_1, B_2, b_7, b_3, B_4, b_5, B_6, B_8$$
$$B_0, b_1, B_2, b_3, B_4, B_6, b_7, b_5, B_8.$$
\item $b_1\DD b_7\DD b_3\DD b_5$. The numberings 
$$B_0, b_1, B_2, B_4, b_5, B_6, b_7, b_3, B_8$$
$$B_0, b_1, B_2, \overline{B_4}, b_5,  b_7, b_3, B_6, B_8$$
(where $\overline{B_4}$ means $B_4$ with order reversed)
have back-edge graphs in which every component is a right-star, and every component is either a clique or left-spike, respectively.
Thus every transversal of the two consists of either a right-star and a clique (and such a transversal has the SRSEH-property
by \ref{cliques}), or a right-star and a left-spike (and such a transversal has the SRSEH-property by \ref{stars}).
Consequently, the result follows from \ref{connected}.
\item $b_1\DD b_7\DD b_5\DD b_3$. Apply \ref{stars} and \ref{connected} to 
$$B_0, b_1, B_2, b_3, B_4, B_6, b_7, b_5, B_8$$
$$B_0, b_1, B_2, b_3, B_4, b_7, b_5, B_6, B_8.$$
\item $b_3\DD b_1\DD b_7\DD b_5$. The numberings
$$B_0, b_3, b_1, \overline{B_2},  B_4, b_5, B_6, b_7, B_8$$
$$B_0, b_1, B_2, b_3, B_4, \overline{B_6}, b_7, b_5, B_8$$
$$B_0, b_3, b_1, B_2,  B_4, b_5, B_6, b_7, B_8$$
have back-edge graphs in which each component is a right-star or clique; each component is a left-star or clique; and each component is a left-star or right-star, respectively.
Hence we may apply \ref{cliques}, \ref{stars} and \ref{connected}.
\item $b_3\DD b_7\DD b_1\DD b_5$. Since $H$ is $D_5$-free, it follows that $B_4$ is null.
Apply \ref{cliques} and \ref{connected} to 
$$B_0,B_2, b_5, b_7, b_1, b_3, B_6, B_8$$
$$B_0,\overline{B_2}, b_5, b_7, b_1, b_3, \overline{B_6}, B_8$$
\end{itemize}

\bigskip

This completes the list of cases with $t=4$; so 
$t\ge 5$, and therefore $t=5$, since $B$ has only three edges and some vertex has degree more than one. The 
subgraph induced on $\{b_1,b_3,b_5,b_7,b_9\}$ has two components, one an edge and the other a three-vertex path.
The three-vertex path might be a monotone path or a left-star or a right-star. Suppose first that it is a monotone path, and so all
components of $B$ are monotone paths. Hence the claim follows if we can exhibit a numbering for which every component of the backedge graph
is a left-star or left-spike, or if there is a numbering for which every component is a right-star or right-spike. 
If the vertices of the three-vertex path in order are $v_i\DD v_j\DD v_k$ where $i<j<k$, and none of $v_{i+1}\LL v_{j-1}$
have positive degree in $B$, then the numbering
$$(v_1\LL v_{i-1}, v_j, v_i\LL v_{j-1}, v_{j+1}\LL v_n)$$
gives a backedge graph in which every component is a left-star, as required; and similarly we may assume that one of 
$v_{j+1}\LL v_{k-1}$ has positive degree in $B$. Consequently the only possibility is (using the $b_i, B_i$ notation)
that the edges of $B$ are $b_1b_5, b_3b_7$ and $b_5b_9$. The backedge graph of the numbering
$$B_0, B_2, b_3, B_4, b_5, B_6, b_7, B_8, b_1,b_9, B_{10}$$
has two components, one an edge and the other a crossed right-star; and the backedge graph of
$$B_0, b_1, b_9, B_2, b_3, B_4, b_5, B_6, b_7, B_8, B_{10}$$
again has two components, one an edge and the other a crossed left-star. Since every component of $B$ itself is a monotone path, the claim follows
from \ref{startri}.

Thus we may assume that a component of $B$ is either a left 2-star or a right 2-star; and from the symmetry under reversal, we may 
assume it is a left 2-star. Let its vertices be $v_i, v_j, v_k$. If none of $v_{i+1}\LL v_{j-1}$ has positive degree in $B$, then 
every component of the backedge graph of the numbering
$$(v_1\LL v_{i-1}, v_{i+1}\LL v_j, v_i, v_{j+1}\LL v_n)$$
is either a right-star or a right-broom, and the claim follows from \ref{stars} and \ref{connected}. 
So we may assume that one of $v_{i+1}\LL v_{j-1}$ has positive degree in $B$. If also some vertex earlier than $v_i$ has positive degree,
the edges of $B$ (in the $B_i,b_i$ notation) are $b_1b_5, b_3b_7, b_3b_9$; and then 
the non-singleton component of the backedge graph of the numbering
$$B_0, b_1, B_2, B_4, b_5, B_6, b_7, B_8, b_9, b_3,B_{10}$$
is a right-bristle, and since every component of $B$ is a left-star with at most three vertices, the claim follows from \ref{bristle}.
So we may assume that $b_1$ is adjacent to exactly two of $b_5,b_7,b_9$, and $b_3$ is adjacent to the third. 
There are three cases, but the same argument applies to each.
Every component of $B$ is a left-star with at most three vertices. The non-singleton component of the backedge graph of the numbering
$$B_0, B_2, b_3, B_4, b_5, B_6, b_7, B_8, b_9, b_1, B_{10}$$
is a crossed right-star; and the non-singleton component of the backedge graph of the numbering
$$B_0, b_1, B_2, b_3, B_4, b_9, \overline{B_8}, b_7, \overline{B_6}, b_5, B_{10}$$
is a left-split, so the claim follows from \ref{split}.
This proves \ref{noD5}.~\bbox

\section{Sparsity}\label{sec:sparsity}
Now we turn to the proof of \ref{mainthm}. We will need the following theorem of R\"odl~\cite{rodl}:
(see for instance~\cite{rodlstrengthen} for this version):
\begin{thm}\label{rodl}
For every graph $H$, and every $\vare>0$, there exists $\delta>0$ such that if $G$ is a graph with no induced subgraph isomorphic
to $H$, there exists $X\subseteq V(G)$ with $|X|\ge \delta |G|$ such that one of the graphs $G[X]$, $\overline{G}[X]$ has maximum degree
less than $\vare\delta|G|$.
\end{thm}
We need a version of this for ordered graphs. To obtain that, we use 
a theorem of R\"odl and Winkler~\cite{RW}, that says:
\begin{thm}\label{orderedornot} For every ordered graph $H$, there exists a graph $P$
such that, for
every numbering of $P$, the resulting ordered graph contains $H$.
\end{thm}

We deduce
\begin{thm}\label{orderedrodl}
For every ordered graph $H$, and every $\vare>0$, there exists $\delta>0$ such that if $G$ is an $H$-free ordered graph,
there exists $X\subseteq V(G)$ with $|X|\ge \delta |G|$ such that one of the graphs $G[X]$, $\overline{G}[X]$ has maximum degree
less than $\vare\delta|G|$.
\end{thm}
\Proof
By \ref{orderedornot}, there is a graph $P$ such that for
every numbering of $P$, the resulting ordered graph contains $H$. Choose $\delta$ as in \ref{rodl}, with $H$ replaced by $P$. We claim that $\delta$
satisfies the theorem. Let $G$ be an ordered graph that does not contain $H$.
From the choice of $P$, it follows that $G$ (as an unordered graph) does not contain
$P$ as an induced subgraph; and so the result follows from the choice of $\delta$. This proves \ref{orderedrodl}.~\bbox

\section{Six-vertex tournaments containing $D_5$}\label{sec:6vertexwithD5}

We will handle the tournaments (with at most six vertices) that contain $D_5$ separately from those that do not, because the
arguments needed are quite different. In this section we handle those that contain $D_5$.
Berger, Choromanski, Chudnovsky and Zerbib~\cite{bergerstrong} proved that $D_5$ itself has the strong EH-property, and we
will show that their proof method also works for what we need.

\begin{figure}[ht]
\centering

\begin{tikzpicture}[scale=1,auto=left]
\def\r{2}
\tikzstyle{every node}=[inner sep=1.5pt, fill=black,circle,draw]
\node (v1) at ({\r*cos(90)}, {\r*sin(90)}) {};
\node (v2) at ({\r*cos(162)}, {\r*sin(162)}) {};
\node (v3) at ({\r*cos(234)}, {\r*sin(234)}) {};
\node (v4) at ({\r*cos(306)}, {\r*sin(306)}) {};
\node (v5) at ({\r*cos(18)}, {\r*sin(18)}) {};
\tikzstyle{every node}=[]
\draw (v1) node [above]           {$1$};
\draw (v2) node [left]           {$2$};
\draw (v3) node [left]           {$3$};
\draw (v4) node [right]           {$4$};
\draw (v5) node [right]           {$5$};

\begin{scope}[thick, decoration={
    markings,
    mark=at position 0.5 with {\arrow{>}}}
    ]
    \draw[postaction={decorate}] (v1) to [bend right=20] (v2);
    \draw[postaction={decorate}] (v2) to [bend right=20] (v3);
    \draw[postaction={decorate}] (v3) to [bend right=20] (v4);
    \draw[postaction={decorate}] (v4) to [bend right=20] (v5);
    \draw[postaction={decorate}] (v5) to [bend right=20] (v1);
    \draw[postaction={decorate}] (v1) to [bend right=20] (v3);
    \draw[postaction={decorate}] (v2) to [bend right=20] (v4);
    \draw[postaction={decorate}] (v3) to [bend right=20] (v5);
    \draw[postaction={decorate}] (v4) to [bend right=20] (v1);
    \draw[postaction={decorate}] (v5) to [bend right=20] (v2);

\end{scope}
\end{tikzpicture}
\caption{$D_5$.} \label{fig:D_5circle}
\end{figure}

Which tournaments do we need to handle? As we said, $D_5$ itself is handled in~\cite{bergerstrong}, so we are concerned with the
tournaments with exactly six vertices that contain $D_5$. With $D_5$ numbered as in figure \ref{fig:D_5circle}, if we add a new
vertex, we can describe
it by giving its set of out-neighbours.
That might be any subset of $\{1\LL 5\}$; but by taking the reverse if necessary, we may
assume the new vertex has at most two out-neighbours (since $D_5$ is isomorphic to its reverse),
and from the symmetry there are only four cases that give nonisomorphic tournaments, namely
$$\emptyset, \{1\}, \{1,2\}, \{1,3\}.$$
The fourth case yields $H_6$, which we cannot do,
so we will just show how to handle the first three. Let us give these three tournaments names: if we start with $D_5$ numbered as above,
and add a new vertex with out-neighbour set $X$ where $X\subseteq \{1\LL 5\}$ (and in-neighbour set $\{1\LL 5\}\setminus X$),
we call the tournament we obtain $D_5^X$. Thus we need to handle $D_5^\emptyset$, $D_5^{\{1\}}$ and $D_5^{\{1,2\}}$.

We need a result proved in~\cite{bergerstrong}. Let us say a digraph $G$ is {\em out-simplicial} if for all distinct $v,x,y\in V(G)$
such that $vx,vy$ are edges, at least one of $xy,yx$ is an edge.
It was proved in~\cite{bergerstrong} that:
\begin{thm}\label{outsimp}
For every out-simplicial digraph $G$ with $|G|>1$, either:
\begin{itemize}
\item there exist disjoint subsets $A,B\subseteq V(G)$ with $|A|,|B|\ge \lfloor |G|/6\rfloor$, such that there are no edges of $G$ between $A,B$ (in either direction); or
\item there exist disjoint subsets $A,B\subseteq V(G)$ with $|A|,|B|\ge \lfloor |G|/6\rfloor$, such that for all $a\in A$ and $b\in B$, there is a directed path in $G$ from $a$ to $b$.
\end{itemize}
\end{thm}

The proof for $D_5$ given in~\cite{bergerstrong} extends to the following:

\begin{thm}\label{boxes}
Let $G$ be a tournament, and let $A_1\LL A_7$ be pairwise disjoint subsets of $V(G)$, each of cardinality at least $W$.
Then either:
\begin{itemize}
\item there exist $1\le i<j\le 7$ such that some vertex in $A_i$ has at least $W/9$ in-neighbours in $A_j$, or some vertex
in $A_j$ has at least $W/9$ out-neighbours in $A_i$; or
\item there is a pure pair in $G$ with order at least $\lfloor W/6\rfloor $; or
\item $G$ contains all of $D_5^\emptyset$, $D_5^{\{1\}}$ and $D_5^{\{1,2\}}$.
\end{itemize}
\end{thm}
\Proof Let $B$ be the graph with vertex set $A_1\cupcup A_7$, where $uv\in E(G)$ if
$u\in A_i$ and $v\in A_j$ for some $i<j$, and $u$ is adjacent from $v$ in $G$. Thus we may assume that
for $1\le i\le 7$, every vertex of $B$ has fewer than $W/9$ $B$-neighbours in $A_i$, for otherwise the theorem holds.
We may assume that $B$ has at least one edge, since otherwise the second outcome holds; so $W\ge 10$.
We begin with the following:
\\
\\
(1) {\em If there exist $v_1,v_2\in A_1$ and $v_4,v_5\in A_6$ such that
$$v_1v_2,v_4v_5, v_5v_1,v_5v_2,v_4v_1,v_2v_4\in E(G),$$
then $G$ contains all of $D_5^\emptyset$, $D_5^{\{1\}}$ and $D_5^{\{1,2\}}$.}
\\
\\
For $1\le i\le 7$ let $A_i'$ be the set of vertices in $A_i$ that are not $B$-adjacent to any of $v_1,v_2,v_4,v_5$. Consequently
$|A_i'|\ge 5W/9$. For each $v_3\in A_3'$, it follows that the subtournament of $G$ induced on $\{v_1\LL v_5\}$ is isomorphic
to $D_5$. (We have chosen the numbering to match that in figure \ref{fig:D_5circle}.) If we choose $v_3\in A_3'$ and $v_6\in A_7'$, not $B$-adjacent (this is possible since $v_3$ has fewer than $W/9$
$B$-neighbours in $A_7'$, and $|A_7'|\ge 5W/9>W/9$), the set
$\{v_1\LL v_6\}$ induces $D_5^\emptyset$. If we choose $v_3\in A_3'$ and $v_6\in A_7'$, $B$-adjacent (this is possible since
otherwise $(A_3', A_7')$ is a pure pair in $G$ and the second outcome of the theorem holds) then
$\{v_1\LL v_6\}$ induces $D_5^{\{1\}}$. If we choose $v_3,v_6\in A_3'$, then $\{v_1\LL v_6\}$ induces $D_5^{\{1,2\}}$. This proves (1).
\\
\\
(2) {\em If there exist $v_1,v_2\in A_1$ and $v_4\in A_3$ and $v_5\in A_6$ such that
$$v_1v_2,v_4v_5, v_5v_1,v_5v_2,v_4v_1,v_2v_4\in E(G),$$
then $G$ contains all of $D_5^\emptyset$, $D_5^{\{1\}}$ and $D_5^{\{1,2\}}$.}
\\
\\
Define $A_1'\LL A_7'$ as before; then again each $|A_i'|\ge 5W/9$. If we choose $v_3\in A_2'$, then $\{v_1\LL v_5\}$
induces $D_5$. If we choose $v_3\in A_2'$ and $v_6\in A_7'$, not $B$-adjacent, then $\{v_1\LL v_6\}$ induces $D_5^\emptyset$.
If we choose  $v_3\in A_2'$ and $v_6\in A_7'$, $B$-adjacent, then $\{v_1\LL v_6\}$ induces $D_5^{\{1\}}$. If we choose
$v_3,v_6\in A_2'$, then $\{v_1\LL v_6\}$ induces $D_5^{\{1,2\}}$. This proves (2).
\\
\\
(3) {\em If there exist $v_1\in A_1$, $v_5\in A_3$ and $v_4\in A_6$, pairwise $B$-adjacent, then
$G$ contains all of $D_5^\emptyset$, $D_5^{\{1\}}$ and $D_5^{\{1,2\}}$.}
\\
\\
For $1\le i\le 7$ let $A_i'$ be the set of vertices in $A_i$ that are not $B$-adjacent to any of $v_1,v_4,v_5$.
Thus each $|A_i'|\ge W-3W/9=2W/3$. If we choose $v_3\in A_2'$ and $v_2\in A_4'$, $B$-adjacent, then $\{v_1\LL v_5\}$
induces $D_5$. If we choose $v_3\in A_2'$ and $v_2\in A_4'$, $B$-adjacent, and choose $v_6\in A_7'$,
not $B$-adjacent to $v_2,v_3$, then $\{v_1\LL v_6\}$ induces $D_5^\emptyset$. If we choose $v_3\in A_2'$ and $v_2\in A_4'$, $B$-adjacent,
and $v_6\in A_5'$, not $B$-adjacent to $v_2,v_3$, then $\{v_1\LL v_6\}$ induces $D_5^{\{1\}}$. Finally, we may assume that some
vertex in $A_4'$ has at least two $B$-neighbours in $A_2'$, because otherwise we can choose $X\subseteq A_2'$ and $Y\subseteq A_4'$
with
$|X|,|Y|= \lfloor W/3\rfloor\ge \lfloor W/6\rfloor$ such that there are no $B$-edges between $X,Y$, and so
$(X,Y)$ is a pure pair of $G$ and the second outcome holds. Let $v_2\in A_4'$ have two $B$-neighbours $v_3,v_6\in A_2'$;
then  $\{v_1\LL v_6\}$ induces $D_5^{\{1,2\}}$. This proves (3).

\bigskip
We assume therefore that none of (1), (2), (3) apply. From now on the argument is exactly as in \cite{bergerstrong}, but we give it for the reader's convenience.
Let $J$ be the digraph with vertex set $A_1$, in which for all distinct $u,v\in A_1$, $v$ is $J$-adjacent from $u$
if $v$ is $G$-adjacent from $u$ and there exists $w\in A_6$ $B$-adjacent to both $u,v$.
\\
\\
(4) {\em $J$ is out-simplicial.}
\\
\\
Suppose that $v_1\in A_1$ is adjacent in $J$ to $v_2,v_2'\in A_1$, and neither of $v_2v_2', v_2'v_2$ is an edge of $J$.
Choose $v_5\in A_6$ $B$-adjacent to $v_1,v_2$, and choose $v_4\in A_6$ $B$-adjacent to $v_1,v_2'$. Since one of $v_1v_2, v_2v_1$
is an edge of $G$, and not an edge of $J$, it follows that $v_4$ is not $B$-adjacent to $v_2$, and $v_5$ is not $B$-adjacent to $v_2'$.
From the symmetry we may assume that $v_4v_5$ is an edge of $G$. But then $v_1,v_2,v_4,v_5$ satisfy the hypotheses of (1),
a contradiction. This proves (4).

From \ref{outsimp}, either
\begin{itemize}
\item there exist disjoint subsets $X,Y\subseteq A_1$ with $|X|,|Y|\ge \lfloor |A_1|/6\rfloor$, such that there are no edges of $J$ between $X,Y$ (in either direction); or
\item there exist disjoint subsets $X,Y\subseteq A_1$ with $|X|,|Y|\ge \lfloor|A_1|/6\rfloor$, such that for all $x\in X$ and $y\in Y$, there is a
directed path in $J$ from $x$ to $y$.
\end{itemize}
Suppose that the first holds. It follows that no vertex in $A_6$ has both a $B$-neighbour in $X$ and a $B$-neighbour in $Y$. Thus either
at least half the vertices in $A_6$  have no $B$-neighbour in $X$, or at least half have no $B$-neighbour in $Y$; and so
$G$ has a pure pair $(P,Q)$ with $P$ one of $X,Y$ and $Q\subseteq A_6$ with $|Q|\ge |A_6|/2\ge W/2$. Since $|P|\ge \lfloor W/6\rfloor$, the second outcome of the theorem holds.

Thus we may assume that the second bullet holds; there exist disjoint subsets $X,Y\subseteq A_1$ with $|X|,|Y|\ge \lfloor |A_1|/6\rfloor$,
such that for all $x\in X$ and $y\in Y$, there is a
directed path in $J$ from $x$ to $y$.
\\
\\
(5) {\em There do not exist $x\in X, y\in Y$ and $z\in A_3$ such that $z$ is $B$-adjacent to $x$,
and $z$ is not $B$-adjacent to $y$.}
\\
\\
Suppose that such $x,y,z$ exist.
Since there is a directed path of $J$ between $x,y$, there is an edge $uv\in E(J)$ such that
$z$ is $B$-adjacent to $u$ and not to $v$. Choose $w\in A_6$ $B$-adjacent to both $u,v$ (this exists from the definition of $J$).
If $z,w$ are not $B$-adjacent, the hypotheses of (2) are satisfied, and if $z,w$ are $B$-adjacent then the hypotheses of (3) are satisfied,
in either case a contradiction.

From (5), either half the vertices in $A_3$ have no $B$-neighbour in $X$, or half the vertices in $A_3$ are $B$-adjacent to all of $Y$;
so there exists $P\subseteq A_3$ with $|P|\ge W/2$ such that one of $(X,P), (P,Y)$ is a pure pair of $G$, and the second outcome of the
theorem holds. This proves \ref{boxes}.~\bbox

We deduce:

\begin{thm}\label{D5mainthm}
If $H$ is a tournament with $|H|\le 6$ that contains $D_5$, and $H$ is different from $H_6,\overline{H_6}$, then $H$
has the strong EH-property.
\end{thm}
\Proof
As we saw, we may assume that $H$ is one of $D_5^\emptyset$, $D_5^{\{1\}}$ and $D_5^{\{1,2\}}$.
Let $J$ be some backedge graph of $H$. By \ref{orderedrodl} there exists $\delta>0$
such that if $G$ is a $J$-free ordered graph,
there exists $X\subseteq V(G)$ with $|X|\ge \delta |G|$ such that one of the graphs $G[X]$, $\overline{G}[X]$ has maximum degree
less than $\delta|G|/126$.
Let $c=\delta/84$; we will show that every $H$-free tournament $G$ with $|G|>1$ has a pure pair of order at least $c|G|$.
Let $G$ be an $H$-free tournament with $|G|>1$. If $|G|\le 1/c$ then $G$ has a pure pair of order 1 that satisfies the theorem, so we
may assume that $|G|>1/c$. Let $(u_1\LL u_m)$ be a numbering of $G$, and let $B$ be its backedge graph.
Thus $B$ is $J$-free.
From the choice of $\delta$, there exists $X\subseteq V(G)$ with $|X|\ge \delta|G|$ such that
one of $B[X]$, $\overline{B}[X]$ has maximum degree
less than $\delta|G|/126$; and by reversing the numbering of $G$ if necessary, we may assume that $B[X]$ has
maximum degree
less than $\delta|G|/126$.

Let $W=6\lceil c|G|\rceil$. Since $c|G|\ge 1$, it follows that
$\lceil c|G| \rceil  \le 2c|G|$, and so
$$\delta|G|/14=6c|G|\le W\le 12c|G|=\delta|G|/7\le |X|/7.$$
Choose
disjoint subsets $A_1\LL A_7$ of $X$, each of cardinality $W$, such that for $1\le i<j\le 7$, if $u_p\in A_i$ and
$u_q\in A_j$ then $p<q$. Since $G$ is $H$-free and $H$ is one of  $D_5^\emptyset$, $D_5^{\{1\}}$ and $D_5^{\{1,2\}}$, it follows
from \ref{boxes} (since $|W|$ is divisible by six) that either
\begin{itemize}
\item there exist $1\le i<j\le 7$ such that some vertex in $A_i$ has at least $W/9$ in-neighbours in $A_j$, or some vertex
in $A_j$ has at least $W/9$ out-neighbours in $A_i$; or
\item there is a pure pair in $G$ with order at least $W/6$.
\end{itemize}
The first is impossible since every vertex in $X$ has degree less than
$\delta|G|/126\le W/9$ in $B[X]$.
Consequently $G$ has a pure pair of order at least $W/6\ge c|G|$. This proves \ref{D5mainthm}.~\bbox


\section{Choosing a backedge graph}\label{sec:6vertexnoD5}

To complete the proof of \ref{mainthm}, we need to handle the six-vertex tournaments that do not contain $D_5$, which is the content
of this section.
We will have to examine all tournaments with at most six vertices, and we will enumerate them by their backedge graphs.
Each tournament may have several different backedge graphs, and we only need to examine one per tournament, so let us try to choose
a good one. In this section we
show that all six-vertex tournaments have backedge graphs with at most four edges; so we can handle most of them by means of 
\ref{noD5}, and the others are handled case by case.
Let us say
a numbering of a tournament $H$ is
{\em optimal} if it has as few backedges as possible, over all numberings of $H$.

\begin{thm}\label{interval}
Let $(v_1\LL v_n)$ be an optimal numbering of a tournament $H$, with backedge graph $B$, and let $1\le i<j\le n$. Then:
\begin{itemize}
\item $v_i$ is $B$-adjacent to at most
$(j-i)/2$ members of $\{v_{i+1}\LL v_j\}$; and $v_j$ is $B$-adjacent to at most
$(j-i)/2$ members of $\{v_{i}\LL v_{j-1}\}$.
\item If $j-i\le 3$ then $B[\{v_i\LL v_j\}]$ has at most one edge.
\item If $j-i=4$ then $B[\{v_i\LL v_j\}]$ has at most three edges. It has three only if they are
$v_iv_{i+4}, v_iv_{i+3}$ and $v_{i+1}v_{i+4}$, and then $G$ contains $D_5$.
\end{itemize}
\end{thm}
\Proof
For the first statement,
$$(v_1\LL v_{i-1},v_{i+1}\LL v_j, v_i, v_{j+1}\LL v_n)$$
is a numbering of $H$, and the
number of its backedges is obtained from the number of backedges of $(v_1\LL v_n)$ by adding the number of $H$-out-neighbours of $v_i$
in $\{v_{i+1}\LL v_j\}$ and subtracting the number of $H$-in-neighbours of $v_i$ in this set; and since $(v_1\LL v_n)$ is optimal,
it follows that at least half of the vertices in $\{v_{i+1}\LL v_j\}$ are $H$-out-neighbours of $v_i$, that is,
$v_i$ is $B$-adjacent to at most
$(j-i)/2$ members of $\{v_{i+1}\LL v_j\}$. This proves half of the first statement and the other half follows from symmetry.

For the second statement, let $j\le i+3$, and suppose that $B[\{v_i\LL v_j\}]$ has at least two edges. Hence there exist
$i\le a<b\le j$ with $b-a<j-i$ such that $v_a, v_b$ are $B$-adjacent; and by the first statement $b-a\ge 2$. Consequently
$j-i=3$. The only pairs of vertices in $\{v_i\LL v_{i+3}\}$ that might be adjacent are $v_iv_{i+2}, v_{i+1}v_{i+3}$ and $v_iv_{i+3}$;
and by the first statement, $v_i$ has at most one $B$-neighbour in $\{v_{i+1}, v_{i+2},v_{i+3}\}$, and $v_{i+3}$
has at most one $B$-neighbour in $\{v_{i}, v_{i+1},v_{i+2}\}$. Thus $v_i$ is not $B$-adjacent to $v_{i+3}$; and so
$v_iv_{i+2}$ and $v_{i+1}v_{i+3}$ are backedges. But then the numbering
$$(v_1\LL v_{i-1}, v_{i+2}, v_{i}, v_{i+3}, v_{i+1}, v_{i+1}\LL v_n)$$
has fewer backedges, a contradiction. This proves the second statement.

For the third, let $j=i+4$, and suppose that $B[\{v_i\LL v_j\}]$ has at least three edges. By the second statement,
$B[\{v_i\LL v_{j-1}\}]$ has at most one edge, and so does $B[\{v_{i+1}\LL v_{j}\}]$; so $B[\{v_i\LL v_j\}]$ has
exactly three edges, and one of them is $v_iv_j$, and each of the other two only appears in one of
$B[\{v_i\LL v_{j-1}\}]$, $B[\{v_{i+1}\LL v_{j}\}]$. Let the other two backedges
be $v_iv_b$ and $v_av_j$ where $a,b\in \{i+1\LL j-1\}$. Now $b\in \{i+2,i+3\}$ and $a\in \{i+1,i+2\}$, from the first statement;
so there are four cases.

\begin{itemize}
\item If $a=b=i+2$, then
$$(v_1\LL v_{i-1}, v_{i+1}, v_{i+4}, v_{i+2}, v_i, v_{i+3}, v_{i+5}\LL v_n)$$
has fewer backedges, a contradiction.
\item If $a=i+1$ and $b=i+2$, then
$$(v_1\LL v_{i-1}, v_{i+2}, v_{i+4}, v_{i}, v_{i+1}, v_{i+3}, v_{i+6}\LL v_n)$$
has fewer backedges, a contradiction; and similarly there is a contradiction if $a=i+2$ and $b=i+3$.
\item If $a=i+1$ and $b=i+3$, then the tournament induced on $\{v_i\LL v_{i+5}\}$ is isomorphic to $D_5$, and
the third outcome of the theorem holds.
\end{itemize}
This proves \ref{interval}.~\bbox

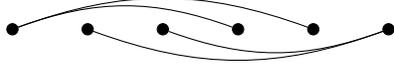
\begin{figure}[H]
\centering

\begin{tikzpicture}[scale=1,auto=left]

\tikzstyle{every node}=[inner sep=1.5pt, fill=black,circle,draw]

\node (v1) at (1,0) {};
\node (v2) at (2,0) {};
\node (v3) at (3,0) {};
\node (v4) at (4,0) {};
\node (v5) at (5,0) {};
\node (v6) at (6,0) {};
\draw (v1) to [bend left=20] (v4);
\draw (v1) to [bend left=20] (v5);
\draw (v2) to [bend right=20] (v6);
\draw (v3) to [bend right=20] (v6);

\end{tikzpicture}

\caption{Backedge graph of $F_6$.} \label{fig:F6}
\end{figure}

Figure \ref{fig:F6} defines the tournament $F_6$.

\begin{thm}\label{minbackedges}
Let $H$ be a tournament with $|H|\le 6$, and let $(v_1\LL v_{|H|})$ be an optimal numbering, with backedge graph $B$.
If $|H|\le 4$, there is at most one
backedge. If $|H|=5$, there are at most three backedges, and at most two unless $H=D_5$. If $|H|=6$, there are at most four backedges,
and at most three unless $H$ is one of $P_7^-, H_6,\overline{H_6}$ and $F_6$.
\end{thm}
\Proof
If $|H|\le 4$ the claim is clear, and if $|H|=5$ the claim follows from \ref{interval}.3 (that is, from the third statement
of \ref{interval}; we will use this notation again).
Thus we may assume that $|H|=6$, and there are at least four backedges.

Suppose that $B[\{v_2\LL v_6\}]$ has at least three edges. By \ref{interval}.3 it has exactly three,
and they are $v_2v_5,v_2v_6$ and $v_3v_6$. All other edges of $B$ are incident with $v_1$. By \ref{interval}.1,
$v_1$ is not $B$-adjacent to $v_6$ or to $v_2$, so the only possible further edges of $B$ are $v_1v_3,v_1v_4$ and $v_1v_5$. By \ref{interval}.1
at most two of them are present; and also by \ref{interval}.1, not both $v_1v_3, v_1v_4$ are backedges. If $v_1v_5$
and one of $v_1v_3,v_1v_4$ are both present then the numbering
$$(v_3,v_4,v_5,v_1,v_6,v_2)$$
has fewer backedges, a contradiction. So exactly one of $v_1v_3,v_1v_4,v_1v_5$ is present. If $v_1v_3$ is a backedge then
$$(v_3,v_1,v_4,v_5,v_6,v_2)$$
has fewer backedges; if $v_1v_4$ is a backedge then
$$(v_4,v_1,v_5,v_6,v_2,v_3)$$
has fewer backedges; and if $v_1v_5$ is a backedge then
$$(v_5,v_1,v_6,v_2,v_3,v_4)$$
has fewer backedges, a contradiction.

Thus we may assume that $B[\{v_2\LL v_6\}]$ has at most two edges. Since by \ref{interval}.1, $v_1$ is incident with at
most two edges of $B$, and $B$ has at least four edges, it follows that exactly two are incident with $b_1$, and $B$ has
four edges altogether. Similarly, exactly two are incident with $v_6$.

Suppose that $v_1,v_6$ are not $B$-adjacent. Since $v_1$ has two $B$-neighbours in $\{v_2\LL v_5\}$, and at most one in $\{v_2,v_3,v_4\}$ by \ref{interval}.1, it follows that $v_1v_5$ is a backedge, and similarly so is $v_2v_6$. Each of $v_1,v_6$
has one further $B$-neighbour; let $v_1v_a$ and $v_bv_6$ be backedges, where $a,b\in \{3,4\}$. Now there are
four cases, $(a,b)=(3,4), (3,3), (4,4), (4,3)$.
\begin{itemize}
\item If $(a,b)=(3,4)$, then
$$(v_3,v_5,v_1,v_6,v_2,v_4)$$
has fewer backedges.
\item If $(a,b)=(3,3)$ then $H$ is isomorphic to $H_6$ (the numbering
$$(v_2,v_3,v_5,v_1,v_4,v_6)$$
gives the backedge graph of figure \ref{fig:H6}). Similarly if $(a,b)=(4,4)$ then $H$ is isomorphic to $\overline{H_6}$.
\item If $(a,b)=(4,3)$ then $H$ is isomorphic to $F_6$ ($B$ itself is the graph of figure \ref{fig:F6}).
\end{itemize}

Thus we may assume that $v_1v_6$ is a backedge. More, we may assume that for {\em every} optimal numbering of $H$, the first
and last vertices are adjacent and both are incident with two backedges for that numbering. Suppose that $v_1v_3\in E(B)$. Then
$$(v_1,v_2,v_3,v_4,v_5,v_6)$$
$$(v_2,v_3,v_1,v_4,v_5,v_6)$$
$$(v_3,v_1,v_2,v_4,v_5,v_6)$$
are all optimal numberings, and so $v_6$ is $B$-adjacent to each of
$v_1,v_2,v_3$, contrary to \ref{interval}.1.

Thus we may assume that $v_1v_3$ is not a backedge, and similarly $v_4v_6$ is not a backedge. Suppose that $v_3v_5$ is a backedge.
Since $v_6$ is incident with two backedges, and \ref{interval}.2 implies that $v_6$ has no $B$-neighbour in $\{v_3,v_4,v_5\}$,
it follows that $v_2v_6$ is a backedge. Also one of $v_1v_4, v_1v_5$ is a backedge. If $v_1v_4$ is a backedge, then $H$ is isomorphic
to $H_6$; and if $v_1v_5$ is a backedge, then the numbering
$$(v_1,v_2,v_4,v_5,v_3,v_6)$$
shows that again $H$ is isomorphic
to $H_6$.

Thus we may assume that $v_3v_5$ is not a backedge, and similarly $v_2v_4$ is not a backedge. But there is a backedge
with both ends in $\{v_2,v_3,v_4,v_5\}$, and so $v_2v_5$ is a backedge. Also $v_1v_b$ is a backedge for some $b\in \{4,5\}$,
and $v_av_6$ is a backedge for some $a\in {2,3}$. There are four cases, $(a,b)=(2,5), (3,5), (2,4), (3,4)$.
\begin{itemize}
\item If $(a,b)=(2,5)$, then $H$ is isomorphic to $F_6$, as we see from the numbering
$$(v_1,v_3,v_4,v_5,v_6,v_2).$$
\item If $(a,b)=(3,5)$, the numbering
$$(v_5,v_1,v_2,v_3,v_4,v_6)$$ is optimal and yet the first and last vertices are not joined by a backedge of this numbering,
a contradiction. Similarly $(a,b)\ne (2,4)$.
\item If $(a,b)=(3,4)$, then $H$ is isomorphic to $P_7^-$.
\end{itemize}
This proves \ref{minbackedges}.~\bbox

Let us observe also that:
\begin{thm}\label{F6}
The tournament $F_6$ has the RSEH-property.
\end{thm}
\Proof
Every component of the backedge graph shown in figure \ref{fig:F6} is a left-star or right-star.
Let $(v_1\LL v_6)$ be the corresponding numbering: then the non-singleton component of the backedge graph of the
numbering $(v_3,v_2,v_4,v_5,v_6,v_1)$ is a clique, so the result follows from \ref{cliques}. This proves \ref{F6}.~\bbox

We deduce our main result \ref{mainthm}, which we restate in a slightly strengthened form:
\begin{thm}\label{mainthm2}
Let $H$ be a tournament with at most six vertices. If $H$ is different from $P_7^-, H_6$ and $\overline{H_6}$ then $H$ has the strong EH-property; and if in addition $H$ is $D_5$-free then $H$ has the rainbow strong EH-property.
\end{thm}
This is immediate from \ref{minbackedges}, \ref{D5mainthm}, \ref{noD5} and \ref{F6}.

\section{Forests and the Paley tournament $P_7$}\label{sec:forests}

We promised earlier to show that $P_7$ does not have the strong EH-property. For that, we use a variant of a theorem of  Erd\H{o}s~\cite{erdos}, which we shall also need in the next section:
\begin{thm}\label{girth}
Let $c>0,g>0$; then there exists  an integer $d>0$ such that for all sufficiently large integers $n$, there is a graph $G$ with $n$ vertices,
such that:
\begin{itemize}
\item every cycle of $G$ has length more than $g$;
\item there do not exist anticomplete $A,B\subseteq V(G)$ with $|A|,|B|\ge cn$; and
\item $G$ has maximum degree less than $d$.
\end{itemize}
\end{thm}
\Proof
Choose an integer $d>0$ with $(dc^2/(8e))^d\ge 6$ (where $e$ is Euler's constant).
Let $n$ be some (sufficiently) large number, and let us take a random graph $G$ with vertex set $\{1\LL 2n\}$, where $i,j$ are adjacent 
independently with probability $p=4/(c^2n)$. Let $x_1$ be the number of pairs $(A,B)$ with $A,B\subseteq V(G)$,
such that $A,B$ are anticomplete, and $|A|,|B|\ge cn$.
Let $x_2$ be the number of cycles in $G$ of length at most $g$; and let $x_3$ be the number of vertices 
with degree at least $d$. We need to estimate the expected value $E(x_i)$ of $x_i$ for $i = 1,2,3$.

First, let $A,B\subseteq V(G)$ be disjoint, with $|A|,|B|\ge cn$.
The probability 
that there are no edges
of $G$ between $A,B$ is at most $(1-p)^{(cn)^2}\le e^{-pc^2n^2}$; and the number of choices of $(A,B)$ is at most $3^{2n}$. So
$$E(x_1)\le e^{-pc^2n^2}3^{2n}\le n/3$$
if $n$ is sufficiently large (since $p=4/(c^2n)$).

The expected number of cycles of length $i$ in $G$ is at most $p^i(2n)^i/(2i)$, so 
$$E(x_2)\le \sum_{3\le i\le g}p^i(2n)^i/(2i)\le p^{g}(2n)^{g}/2\le (8/c^2)^g/2\le n/3.$$ 

For a vertex $v$, the probability that $v$ has degree at least $d$ is at most 
$ \binom{2n}{d}p^d\le (2pn)^d/d!$; and since $d!\ge  (d/e)^d$ by Stirling's formula, it follows that
 the probability that $v$ has degree at least $d$ is at most 
$(8e/(c^2d))^d\le 1/6$. So $E(x_3)\le n/3$.

Hence the expected value of $x_1+x_2+x_3$
is at most $n$; and so there is a choice of $G$
where $x_1+x_2+x_3\le n$.
Hence by deleting $n$ vertices appropriately we obtain a graph with $n$ vertices as in the theorem.
This proves \ref{girth}.~\bbox

Now we can prove \ref{forest}, which we restate:
\begin{thm}\label{forest2}
Let $H$ be a tournament with the strong EH-property. Then there is a numbering of $H$ such that the backedge graph is a forest.
Consequently there is a partition of $V(H)$ into two subsets both inducing transitive tournaments.
\end{thm}
\Proof
Choose $c>0$ such that every $H$-free tournament $G$ with $|G|>1$ admits a pure pair with order at least $c|G|$.
Choose $d$ satisfying \ref{girth} with $c$ replaced by $c/2$ and $g$ replaced by $|H|$.
Let $n\ge 2d/c$ be some large number, large enough that there is a graph $J$ with $n$ vertices, satisfying the three bullets of \ref{girth}
with $c$ replaced by $c/2$ and $g$ replaced by $|H|$.
Take a numbering $(v_1\LL v_n)$ of $J$, and let $G$ be the tournament such that $J$ is the backedge graph of $G$ under 
this numbering.
Suppose that there is a pure pair $(X,Y)$ in $G$ with $|X|,|Y|\ge cn$. Choose $i$ minimum such that $|\{v_1\LL v_i\}\cap X|\ge cn/2$,
and let $A=\{v_1\LL v_i\}\cap X$ and $B=Y\cap \{v_{i+1}\LL v_n\}$. Since $J$ has maximum degree less than $d\le cn/2$,
and $v_{i}\in X$ is $J$-adjacent to every vertex in $Y\setminus B$, it follows that $|B|\ge cn/2$; and yet $A,B$ are anticomplete,
contrary to the choice of $J$.

Thus $G$ has no pure pair 
$(X,Y)$ in $G$ with $|X|,|Y|\ge cn$. From the definition of $c$, it follows that $G$ contains $H$. The backedge graph
for $H$ under the numbering induced by $(v_1\LL v_n)$ has no cycles, since all cycles of $J$ have length more than $|H|$.
Hence it is a forest. 

This forest is two-colourable; and each colour class induces a transitive subtournament of $H$, since it is a stable set of a
backedge graph of $H$. This proves \ref{forest}.~\bbox

We deduce: 
\begin{thm}\label{badP_7}
$P_7$ does not have the strong EH-property.
\end{thm}
\Proof
It suffices to show that the vertex set of $P_7$ cannot be partitioned as in \ref{forest}, and to show that it suffices
to show that $P_7$ has no four-vertex transitive subtournament. But for every vertex of $P_7$, its three out-neighbours form
a cyclic triangle. This proves \ref{badP_7}.~\bbox

In \cite{pure1} we mentioned \ref{forest} and several other conditions that were necessary if a tournament is to have the strong EH-property.
But we subsequently observed that each of the other conditions was implied by the first; and at the moment,
\ref{forest} is the only necessary condition we know. As we mentioned in the introduction, it might be that having a backedge graph that is a forest is necessary and 
sufficient for a tournament to have the strong EH-property. One
piece of evidence in favour of this is the following, which follows from results of \cite{pure6}:
\begin{thm}\label{pure6}
For a tournament $H$, the following are equivalent:
\begin{itemize}
\item some backedge graph of $H$ is a forest;
\item for every $c>0$ there exists $\vare>0$ such that for every $H$-free tournament $G$ with $|G|>1$, 
there is a pure pair in $G$ of order at least 
$\vare|G|^{1-c}$.
\end{itemize}
\end{thm}
\Proof
It is shown in \cite{pure6} that:
\\
\\
(1) {\em If $J$ is an ordered forest, then 
for all $c>0$, there exists $\vare>0$ such that
if $G$ is an ordered graph with $|G|>1$ that is both $J$-free and $\overline{J}$-free, then
$G$ contains a pure pair of order at least $\vare |G|^{1-c}$.
}
\\
\\
Now let $H$ be a tournament. Suppose first that some backedge graph $J$ of $H$ is a forest. Then $\overline{J}$ is also a backedge graph of $H$
(reversing the numbering); and if $G$ is an $H$-free tournament, and $B$ is its backedge graph, then
$B$ contains neither $J$ nor $\overline{J}$, and so (1) implies that $B$ has the desired pure pair, and hence, by \ref{purepair}, 
so does $G$.

For the converse, let $H$ be a tournament for which no backedge graph is a forest. Let $c<1/|H|$; we claim there is no $\vare$ satisfying
the second bullet of the theorem. Let $\vare>0$.
An argument like that of \ref{girth} shows that if we take a
random graph $J$ on $n$ vertices where $n$ is sufficiently large, in which every edge is present independently with probability
$\frac12 n^{-1+1/|H|}$, then with high probability, there will be a set $X$ of at least $n/2$ vertices in which $J[X]$ has no cycle
of length at most $|H|$ and has no pure pair of order at least $\vare|X|^{1-c}/2$. Number $X$ arbitrarily, and let $G$ be the tournament
with $J$ (and this numbering) as a backedge graph. Then $G$ does not contain $H$, since if it did, the induced numbering of $H$
would have backedge graph contained in $J$ with a cycle of length at most $|H|$. And yet $G$ has no pure pair of order at least
$\vare |G|^{1-c}$, by \ref{purepair}, and so $\vare$ does not satisfy the second bullet of the theorem. This proves \ref{pure6}.~\bbox

\section{$D_5$ and $P_7^-$ do not have the rainbow strong EH-property}\label{sec:P7bad}

We claimed earlier that $D_5$ does not have the RSEH-property. The same holds for $P_7^-$, and even excluding them both 
simultaneously is not enough. We will show:
\begin{thm}\label{nonRSEH}
For all $c>0$, and infinitely many integers $n$, 
there is a tournament $G$ with $n$ vertices, and a blockade $\mathcal{B}$ in $G$ of length at least $1/c$, such that
$G$ has no pure pair of order at least $cW(\mathcal{B})$, and contains no $\mathcal{B}$-rainbow copy of either of $D_5, P_7^-$.
\end{thm}

To show this we need a construction as follows. Let $G$ be an ordered graph. A {\em walk} in $G$ is a sequence
$$p_0,p_1\LL , p_r,$$ 
where $p_0\LL p_r\in V(G)$ and there is an edge of $G$ with ends $p_{i-1}, p_i$ for $1\le i\le r$. (We do not require $p_0\LL p_r$
all to be distinct, but consecutive terms are distinct.) Its {\em length} is $r$, and its {\em imbalance}  is $N_1-N_2$, where
$N_1$ is the number of $i\in \{1\LL r\}$ such that $p_{i-1}$ is before $p_{i}$ in the numbering of $G$, and $N_2$
is the number of $i$ such that $p_{i-1}$ is after $p_{i}$. 
A walk is {\em balanced} if its imbalance is zero, and {\em unbalanced} otherwise; and {\em closed} if $p_0=p_r$.

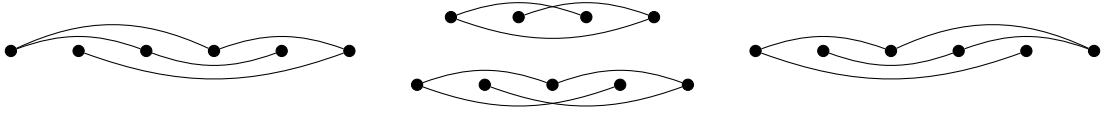
\begin{figure}[H]
\centering

\begin{tikzpicture}[scale=.9,auto=left]

\tikzstyle{every node}=[inner sep=1.5pt, fill=black,circle,draw]

\node (w1) at (1.5,1/2) {};
\node (w2) at (2.5,1/2) {};
\node (w3) at (3.5,1/2) {};
\node (w4) at (4.5,1/2) {};
\draw (w1) to [bend left=20] (w3);
\draw (w2) to [bend left=20] (w4);
\draw (w1) to [bend right=20] (w4);

\node (v1) at (1,-1/2) {};
\node (v2) at (2,-1/2) {};
\node (v3) at (3,-1/2) {};
\node (v4) at (4,-1/2) {};
\node (v5) at (5,-1/2) {};
\draw (v1) to [bend left=20] (v3);
\draw (v3) to [bend left=20] (v5);
\draw (v1) to [bend right=20] (v4);
\draw (v2) to [bend right=20] (v5);

\node (u1) at (6,0) {};
\node (u2) at (7,0) {};
\node (u3) at (8,0) {};
\node (u4) at (9,0) {};
\node (u5) at (10,0) {};
\node (u6) at (11,0) {};
\draw (u1) to [bend left=20] (u3);
\draw (u1) to [bend right=20] (u5);
\draw (u3) to [bend left=25] (u6);
\draw (u2) to [bend right=20] (u4);
\draw (u4) to [bend left=20] (u6);

\node (u1) at (0,0) {};
\node (u2) at (-1,0) {};
\node (u3) at (-2,0) {};
\node (u4) at (-3,0) {};
\node (u5) at (-4,0) {};
\node (u6) at (-5,0) {};
\draw (u1) to [bend right=20] (u3);
\draw (u1) to [bend left=20] (u5);
\draw (u3) to [bend right=25] (u6);
\draw (u2) to [bend left=20] (u4);
\draw (u4) to [bend right=20] (u6);

\end{tikzpicture}

\caption{Ordered graphs for \ref{random}.} \label{fig:oddpaths}
\end{figure}

\begin{thm}\label{random}
Let $k\ge 1$ be an integer, and let $c>0$. Then there is an integer $D$, such that for all sufficiently large integers $W$, 
there is an ordered graph $J$ with $kW$ vertices, and the following properties:
\begin{itemize}
\item every vertex has degree at most $D$;
\item $G$ admits a respectful blockade $\mathcal{B}=(B_1\LL B_k)$ of width $W$;
\item $G$ has no pure pair of order at least $cW$;
\item every closed walk in $J$ of length at most six is balanced;
\item there is no $\mathcal{B}$-rainbow copy in $J$ of any of the ordered graphs shown in figure \ref{fig:oddpaths}.
\end{itemize}
\end{thm}

\Proof 
Let $c'=c/k$, and $g=6\cdot 3^k$.
Choose $d$ to satisfy \ref{girth} with $c$ replaced by $c'$. Let $D=d^{3^k}$. Let $W$ be a sufficiently large integer. Then by \ref{girth} there is a 
graph $J_k$
with $kW$ vertices $v_1\LL v_{kW}$, such that
\begin{itemize}
\item every cycle of $J_k$ has length more than $g$;
\item there is no anticomplete pair in $J_k$ of order at least $c'kW=cW$; and
\item $J_k$ has maximum degree less than $d$.
\end{itemize}

For $1\le i\le k$ let $B_i=\{v_j:(i-1)W<j\le iW\}$, and $\mathcal{B}=(B_1\LL B_k)$. If $u,v\in V(J_k)$, and $u\in B_i$, $v\in B_j$, we define 
the {\em $\mathcal{B}$-length} of the pair $(u,v)$ to 
be $|j-i|$, and the {$\mathcal{B}$-length} of an edge $uv$ is the $\mathcal{B}$-length of $(u,v)$. We say $P$ is a {\em welcoming path}
if 
\begin{itemize}
\item $P$ is a path of length three with $V(P)\subseteq V(J_k)$, 
with ends $s,t$ where $s$ is before $t$ in the numbering $(v_1\LL v_{kW})$, 
\item the $\mathcal{B}$-length of $(s,t)$ is at least one;
\item every edge of $P$ has $\mathcal{B}$-length strictly greater than the $\mathcal{B}$-length of $(s,t)$; and
\item the walk of length three from $s$ to $t$ in $P$ has imbalance one (note that this is different from having imbalance $-1$).
\end{itemize}

For $i=k-1\LL 1$ we define $J_i$ as follows. We say a pair $(s,t)$ of vertices of $J_{i+1}$ is {\em $i$-good} if $s$ is before $t$ 
in the numbering $(v_1\LL v_{kW})$, and the $\mathcal{B}$-length of $(s,t)$ is exactly $i$, and $s,t$ are nonadjacent in $J_{i+1}$, 
and there is a welcoming path in $J_{i+1}$ (not necessarily induced) with ends $s,t$.
We construct $J_{i}$ from $J_{i+1}$ by adding an edge between $s,t$
for every $i$-good pair $(s,t)$.

Let $J=J_1$; we claim that $J$ satisfies the theorem. First, let $d_i$ denote the maximum degree
of $J_i$; then since each vertex of $J_i$ is an end vertex of at most $d_{i+1}(d_{i+1}-1)^2$ paths of length three, it follows
that $d_{i}\le d_{i+1}(d_{i+1}-1)^2+d_{i+1}\le d_{i+1}^3$. Since $d_k\le d$, it follows that $J$ has maximum degree at most $d^{3^k}=D$.

Second, since $J_k$ has no anticomplete pair of order at least $c'kW$, the same holds for $J$. Third, for every closed walk
of $J_{i}$, we can replace each edge $e$ of $J_{i}$ not in $J_{i+1}$ by a three-edge walk along the corresponding welcoming path $P$
of $J_{i+1}$; and since this three-edge walk has imbalance the same as the corresponding one-edge walk along $e$, it follows that
there is a closed walk of $J_{i+1}$ with the same imbalance and with length at most three times as great. Since
every cycle of $J_k$ has length more than $g$, and so every closed walk in $J_k$ with length at most $g$ is balanced,
it follows that 
every closed walk of $J_i$ of length at most $g3^{i-k}$ is balanced, and in particular every closed walk of $J$ with length at most
$g3^{-k}=6$ is balanced.

Fourth, we must show that $J$ contains no $\mathcal{B}$-rainbow copy of any of the four graphs in figure \ref{fig:oddpaths}.
For this we use:
\\
\\
(1) {\em For every welcoming path of $J$, its ends are adjacent in $J$.}
\\
\\
Let $P$ be a welcoming path of $J$, with ends $s,t$ where $s$ is earlier than $t$. Let $i$ be the $\mathcal{B}$-length of $(s,t)$.
Since every edge of $P$ has $\mathcal{B}$-length more than $i$, it follows that every such edge is an edge of $J_{i+1}$ (because
all edges added later have $\mathcal{B}$-length at most $i$), and so $P$ is a welcoming path of $J_{i+1}$. But then $s,t$ are 
adjacent in $J_i$ and hence in $J$. This proves (1).

\bigskip

Suppose that $J$ contains a $\mathcal{B}$-rainbow copy $H$ of one of the four graphs in figure \ref{fig:oddpaths}. (Thus $H$ is induced.)
Suppose first that $|H|=4$, and let its numbering be $(u_1,u_2,u_3,u_4)$. Thus its edges are $u_1u_3,u_2u_4, u_1u_4$.
Let $i$ be the $\mathcal{B}$-length of $(u_2,u_3)$; thus $i\ge 1$ since $H$ is $\mathcal{B}$-rainbow, and for the same reason, 
all three edges of $H$ have $\mathcal{B}$-length at least $i+1$. But then $H$ is a welcoming path of $J$ and its ends are nonadjacent,
contrary to (1).

Now suppose that $|H|=5$, and so its edges are $u_2\DD u_5\DD u_3\DD u_1\DD u_4$, where $(u_1\LL u_5)$ is its numbering. Let the $\mathcal{B}$-length of $(u_1,u_2)$
be $i_1$, and that of $(u_4,u_5)$ be $i_2$. From the symmetry we may assume that $i_1\le i_2$. But then all edges of the path
$u_1\DD u_3\DD u_5\DD u_2$ have $\mathcal{B}$-length strictly more than $i_1$, since $H$ is $\mathcal{B}$-rainbow; so it is welcoming,
contrary to (1).

Finally, suppose that $|H|=6$; and from the symmetry, we may assume that its edges are $u_2\DD u_4\DD u_6\DD u_3\DD u_1\DD u_5$, where
$(u_1\LL u_6)$ is its numbering. Let the $\mathcal{B}$-length of $(u_2,u_3)$
be $i_1$, and that of $(u_5,u_6)$ be $i_2$. If $i_1\le i_2$, then the path $u_2\DD u_4\DD u_6\DD u_3$ is welcoming, and if
$i_2\le i_1$ then the path $u_5\DD u_1\DD u_3\DD u_6$ is welcoming, and in either case we have a contradiction to (1).
This proves \ref{random}.~\bbox

We deduce \ref{nonRSEH2}, which we restate:
\begin{thm}\label{nonRSEH2}
For all $c>0$, and infinitely many integers $n$,
there is a tournament $G$ with $n$ vertices, and a blockade $\mathcal{B}$ in $G$ of length at least $1/c$, such that
$G$ has no pure pair of order at least $cW(\mathcal{B})$, and there is no $\mathcal{B}$-rainbow copy of either of $D_5, P_7^-$.
Conequently $D_5, P_7^-$ do not have the RSEH-property.
\end{thm}
\Proof Let $k=\lceil 2/c\rceil$, and choose $D$, and $W>2D/c$ sufficiently large that the construction $J$ of \ref{random} exists 
with $c$ replaced by $c/2$. 
Let $G$ be the tournament with backedge graph $J$. Since every vertex of $J$ has degree at most $D<cW/2$, it follows that
$J$ has no pure pair of order at least $cW/2$, and so $G$ has no pure pair of order at least $cW$, by \ref{purepair}. 
By examining all the backedge graphs of $D_5$ (there are 24 of them)
and all the backedge graphs of $P_7^-$ (there are 240 of them)
we observe that each of them contains an unbalanced cycle of length at most five, or one of the ordered graphs of figure
\ref{fig:oddpaths}. Consequently there is no $\mathcal{B}$-rainbow copy in $J$ of any backedge graph of $D_5$ or of $P_7^-$,
and so $G$ contains no $\mathcal{B}$-rainbow copy of $D_5$ or of $P_7^-$. This proves \ref{nonRSEH2}.~\bbox


\begin{thebibliography}{99}
\bibitem{aps} N. Alon, J. Pach and J. Solymosi, ``Ramsey-type theorems with forbidden subgraphs'', {\em Combinatorica} {\bf 21} (2001), 155--170.

\bibitem{heroes} E. Berger, K. Choromanski, M. Chudnovsky, J. Fox, M. Loebl, A.
Scott, P. Seymour and S. Thomass\'e, ``Tournaments and colouring'',
{\em J. Combinatorial Theory, Ser. B}, {\bf 103} (2013), 1--20.


\bibitem{galaxies} E. Berger, K. Choromanski and M. Chudnovsky, ``Forcing large transitive subtournaments'', 
{\em J. Combinatorial Theory, Ser. B}, {\bf 112} (2015), 1--17.

\bibitem{berger} E. Berger, K. Choromanski and M. Chudnovsky, ``On the Erd\H{o}s-Hajnal conjecture for six-vertex tournaments'',
{\em European Journal of Combinatorics}, {\bf 75}
(2019) 113--122,
{\tt arXiv:1508.04992}.

\bibitem{bergerstrong} E. Berger, K. Choromanski, M. Chudnovsky and S. Zerbib, ``Tournaments and the strong Erd\H{o}s-Hajnal 
property'', {\em European J. Combinatorics}, {\bf 100} (2022), 103440, {\tt arXiv:2002.07248}.


\bibitem{constellations} K. Choromanski, ``EH-suprema of tournaments with no nontrivial homogeneous sets'', {\em J. Combinatorial Theory,
Ser. B}, {\bf 114} (2015), 97--123.

\bibitem{EHC5} M. Chudnovsky, A. Scott, P. Seymour and S. Spirkl, ``Erd\H{o}s-Hajnal for graphs with no 5-hole'', 
{\em Proceedings of the London Math. Soc.}, {\bf 126} (2023), 997--1014,
{\tt arXiv:2102.04994}.

\bibitem{pure1} M. Chudnovsky, A. Scott, P. Seymour and S. Spirkl,
``Pure pairs. I. Trees and linear anticomplete pairs'',
{\em Advances in Math.}, {\bf 375} (2020), 107396, {\tt arXiv:1809.00919}. 

\bibitem{rodlstrengthen} M. Chudnovsky, A. Scott, P. Seymour and S. Spirkl, ``Strengthening R\"odl's theorem'', submitted for publication, {\tt arXiv:2105.07370}.

\bibitem{erdos} P. Erd\H{o}s, ``Graph theory and probability'', {\em Canad. J. Math.} {\bf 11} (1959), 34--38.

\bibitem{EH0} P. Erd\H{o}s and A. Hajnal, ``On spanned subgraphs of graphs'',
{\em Contributions to Graph Theory and its
Applications} (Internat. Colloq., Oberhof, 1977) (German), 80--96, Tech.
Hochschule Ilmenau, Ilmenau, 1977.

\bibitem{EH}  P. Erd\H{o}s and A. Hajnal, ``Ramsey-type theorems'',
{\em  Discrete Applied Math.} {\bf 25} (1989), 37--52.

\bibitem{fox} J. Fox, ``A bipartite analogue of Dilworth's theorem'', {\em Order} {\bf 23} (2006), 197--209.

\bibitem{density4} T. Nguyen, A. Scott and P. Seymour, ``Induced subgraph density. IV. New graphs with the Erd\H{o}s-Hajnal
property'',
manuscript May 2023.

\bibitem{rodl} V. R\"odl, ``On universality of graphs with uniformly distributed edges'',
{\em Discrete Math.} {\bf 59} (1986), 125--134.

\bibitem{RW} V. R\"odl and P. Winkler, ``A Ramsey-type theorem for orderings of a graph'', {\em SIAM Journal
of Discrete Mathematics} {\bf 2} (1989), 402--406.

\bibitem{pure5}  A. Scott, P. Seymour and S. Spirkl, ``Pure pairs. V. Excluding some long subdivision'', 
{\em Combinatorica}, to appear, {\tt arXiv:2105.03956}.


\bibitem{pure6}  A. Scott, P. Seymour and S. Spirkl, ``Pure pairs. VI. Excluding an ordered tree'',  {\em SIAM J Disc Math.}, {\bf 36} (2022), 170--187, {\tt arXiv:2009.10671}.

\bibitem{7vertex} S. Zayat and S. Ghazal, ``About the Erd\H{o}s-Hajnal conjecture for seven-vertex
tournaments'', {\tt arXiv:2010.12331}.

\end{thebibliography}
\end{document}